\documentclass{amsart}
\usepackage{amsmath,amssymb,amsfonts,amsthm}

\newcommand{\setC}{\mathbb{C}}
\newcommand{\setN}{\mathbb{N}}

\newcommand{\setZ}{\mathbb{Z}}

\newcommand{\inner}[3][]{( #2|#3 )_{#1}}   % prodotto scalare
\newcommand{\del}{\partial}

\newcommand{\ud}{\mathrm{d}}

\newcommand{\tp}[1]{\sp{\otimes #1}} % tensor power
\DeclareMathOperator{\Aut}{Aut} % Automorphisms

\DeclareMathOperator{\diag}{diag} %diagonal

\DeclareMathOperator{\tr}{tr}

\DeclareMathOperator{\Coeff}{Coeff} % coefficiente di Laurent

\newcommand{\onehalf}{\frac{1}{2}}

\newcommand{\U}{\sqrt{-1}}  % radice di -1
\newcommand{\R}{\mathcal{R}}  % R calligafica

\newcommand{\abs}[1]{\lvert#1\rvert}
\newcommand{\norm}[1]{\lVert#1\rVert}
\newcommand{\card}[1]{\lvert#1\rvert}

\let\geq=\geqslant % come dice Campanato, il segno di '>' sottolineato {\'e}
\let\leq=\leqslant % stato introdotto solo perch{\`e} i tipografi non
\newcommand{\moduli}{\mathcal{M}}           % spazio dei moduli di curve
\newcommand{\modulibar}{\overline{\moduli}} % compattificazione
\theoremstyle{plain}

\newtheorem*{theorem*}{Theorem}
\newtheorem*{mainthm*}{Main Theorem}

\newtheorem*{lemma*}{Lemma}

\theoremstyle{remark}

\theoremstyle{definition}

\numberwithin{equation}{section} % numerazione delle equazioni
\usepackage[arc,all,knot,frame,curve,poly]{xy}

\newcommand{\savexyscale}{\let\saved@everyxy=\everyxy}
\newcommand{\restorexyscale}{\let\everyxy=\saved@everyxy}
\newcommand{\xyscale}[1]{\everyxy={\POS#1,0,}}

\xyscale{/r7mm/:}

\providecommand{\rgholescalei}{1.4}
\providecommand{\rgholescaleii}{1.6}    % => 2.24
\providecommand{\rgholescaleiii}{1.9}   % => 4.25

\providecommand{\rgvertexscalei}{0.75}
\providecommand{\rgvertexscaleii}{1.3}  % => 0.97
\providecommand{\rgvertexmark}{\ensuremath\bullet}

\newcommand{\rghole}[2][]{% 
\POS{0*+{#1}="CENTER"},% 
\ifnum #2=-1% 
\rgholeoneup% 
\else\ifnum #2=1% 
\rgholeone% 
\else\ifnum #2=2% 
\rgholetwo% 
\else% 
\rgholepoly #2% 
\fi\fi\fi% 
} 

\newcommand{\rgholeone}{% 
        \POS{"CENTER"*\cir(1,0){}}%
        \POS{DC*{\rgvertexmark}}="VERTEX1";%
        \POS{"CENTER";"VERTEX1"**\dir{}?(\rgholescalei)}="AUXi1";%
        \POS{"CENTER";"AUXi1"**\dir{}?(\rgholescaleii)}="AUXii1";%
        \POS{"CENTER";"AUXii1"**\dir{}?(\rgholescaleiii)}="AUXiii1";%
}

\newcommand{\rgholeoneup}{% 
        \POS{"CENTER"*\cir(1,0){}}%
        \POS{UC*{\rgvertexmark}}="VERTEX1";%
        \POS{"CENTER";"VERTEX1"**\dir{}?(\rgholescalei)}="AUXi1";%
        \POS{"CENTER";"AUXi1"**\dir{}?(\rgholescaleii)}="AUXii1";%
        \POS{"CENTER";"AUXii1"**\dir{}?(\rgholescaleiii)}="AUXiii1";%
}

\newcommand{\rgholetwo}{% 
        \POS{"CENTER"*\cir(1,0){}}="@o",%
        \POS{"@o",UC*{\rgvertexmark}}="VERTEX1",%
        \POS{"@o",DC*{\rgvertexmark}}="VERTEX2",%
        \POS{"CENTER";"VERTEX1"**\dir{}?(\rgholescalei)}="AUXi1";%
        \POS{"CENTER";"AUXi1"**\dir{}?(\rgholescaleii)}="AUXii1";%
        \POS{"CENTER";"AUXii1"**\dir{}?(\rgholescaleiii)}="AUXiii1";%
        \POS{"CENTER";"VERTEX2"**\dir{}?(\rgholescalei)}="AUXi2";%
        \POS{"CENTER";"AUXi2"**\dir{}?(\rgholescaleii)}="AUXii2";%
        \POS{"CENTER";"AUXii2"**\dir{}?(\rgholescaleiii)}="AUXiii2";%
}

\newcommand{\rgholepoly}[1]{% 
        \POS{\xypolygon#1"VERTEX"{% 
        ~>{-}% linee singole per i lati 
        \rgvertexmark 
        }}="CENTER"; 
        \POS{\xypolygon#1"AUXi"{% 
        ~:{(\rgholescalei,0):}% ingrandisce 
        ~>{}% lati vuoti 
        }}; 
        \POS{\xypolygon#1"AUXii"{% 
        ~:{(\rgholescaleii,0):}% ingrandisce 
        ~>{}% lati vuoti 
        }}; 
}

\newcommand{\rgvertex}[2][\relax]{%
\ifx#1\relax%
        \POS*{\displaystyle\rgvertexmark}="CENTER",%
\else%
        \POS*+{\displaystyle #1}*\cir{}="CENTER",% DECOR nel centro
\fi%
\ifnum #2=-1% 
\rgvertexoneup% 
\else\ifnum #2=1% 
\rgvertexone% 
\else\ifnum #2=2% 
\rgvertextwo% 
\else% 
\rgvertexpoly{#2}% 
\fi\fi\fi% 
} 

\newcommand{\rgvertexone}{% 
\POS{DC}="VERTEX1";%
\POS{"CENTER";"VERTEX1":(\rgvertexscalei,0):(1,0)}="AUXi1";%
\POS{"CENTER";"AUXi1":(\rgvertexscaleii,0):(1,0)}="AUXii1";%
}

\newcommand{\rgvertexoneup}{% 
\POS{UC}="VERTEX1";%
\POS{"CENTER";"VERTEX1"**\dir{}?(\rgvertexscalei)}="AUXi1";%
\POS{"CENTER";"AUXi1"**\dir{}?(\rgvertexscaleii)}="AUXii1";%
}

\newcommand{\rgvertextwo}{% 
\POS{DC}="VERTEX1";%
\POS{"CENTER";"VERTEX1"**\dir{}?(\rgvertexscalei)}="AUXi1";%
\POS{"CENTER";"AUXi1"**\dir{}?(\rgvertexscaleii)}="AUXii1";%
\POS{"VERTEX1";"CENTER"**\dir{}?(2.0)}="VERTEX2",%
\POS{"CENTER";"VERTEX2"**\dir{}?(\rgvertexscalei)}="AUXi2";%
\POS{"CENTER";"AUXi2"**\dir{}?(\rgvertexscaleii)}="AUXii2";%
}

\newcommand{\rgvertexpoly}[1]{\POS{% 
\xypolygon#1"VERTEX"{~:{R:}~<{-}~>{}}%,="CENTER"; 
\xypolygon#1"AUXi"{%
~:{(\rgvertexscalei,0):}% ingrandisce
~>{}% lati vuoti 
}%;
\xypolygon#1"AUXii"{% 
~:{(\rgvertexscaleii,0):}% ingrandisce 
~>{}% lati vuoti 
}%, 
}}

\newcommand{\missing}[2][\cdots]{\POS{%
        0;"AUXi#2"**\dir{}?(.9)*{\displaystyle{#1}}%
}}

\newcommand{\loose}[1]{\POS{% 
"VERTEX#1";"AUXii#1"**\dir{-},% 
}}

\newcommand{\fence}[3]{\POS{% 
"VERTEX#1"\PATH~={**\dir{-}} '"AUXii#1"%
^>{#3}%
_>{#2}
}}

\newcommand{\cilia}[3][\bigstar]{\POS{%
    "AUXi#2";"AUXi#3"**\dir{}?(.5)="m@",%
    "CENTER";"m@"**\dir{}?(1.0)*{#1}%
}}

\newcommand{\cilialontano}[3][\bigstar]{\POS{%
    "AUXi#2";"AUXi#3"**\dir{}?(.5)="m@",%
    "CENTER";"m@"**\dir{}?(1.8)*{#1}%
}}

\newcommand{\ciliamedio}[3][\bigstar]{\POS{%
    "AUXi#2";"AUXi#3"**\dir{}?(.5)="m@",%
    "CENTER";"m@"**\dir{}?(1.5)*{#1}%
}}

\newcommand{\join}[3][]{\POS{% 
"AUXi#3";"AUXi#2"**\dir{}?="m@",% 
"CENTER";"m@"**\dir{}?(2.0)="c@",% 
"VERTEX#3";"VERTEX#2"**\crv{~**\dir{-}% 
  "AUXi#3"&"AUXii#3"&"c@"&"AUXii#2"&"AUXi#2"}%
  ?*{#1}%
}} 

\newcommand{\genericfences}[1][n]{%
\fence{1}{\ldots}{i_{#1}}%
\fence{2}{i_{#1}}{i_1}%
\fence{3}{i_1}{i_2}%
\fence{4}{i_2}{\ldots}%
\missing[\cdots]{5}
}

\def\xyc#1\endxyc{{\xy*!C\xybox{#1}\endxy}}
\providecommand{\Xy}{\leavevmode%
 \hbox{\kern-.1em X\kern-.3em\lower.4ex\hbox{Y\kern-.15em}}}

\usepackage{prettyref}
\newrefformat{appendix}{Appendix~\ref{#1}}
\newrefformat{prop}{Proposition~\ref{#1}}
\newrefformat{theorem}{Theorem~\ref{#1}}
\newrefformat{lemma}{Lemma~\ref{#1}}
\newrefformat{xmp}{Example~\ref{#1}}
\newrefformat{definition}{Definition~\ref{#1}}
\newrefformat{cor}{Corollary~\ref{#1}}
\newrefformat{rem}{Remark~\ref{#1}}

\newcommand{\SN}[1][\relax]{\mathcal{S}%
  \ifx#1\relax\relax\else\sb{#1}\fi}
\newcommand{\M}{\moduli}
\newcommand{\Mbar}{\modulibar}
\newcommand{\Hermitian}[1][\relax]{{\mathcal
H}\ifx#1\relax\relax\else(#1)\fi}

\newcommand{\Vertices}[1]{#1^{(0)}}
\newcommand{\Edges}[1]{#1^{(1)}}
\newcommand{\Holes}[1]{#1^{(2)}}

\newcommand{\gint}[2][\Lambda]{%
\Biggl\langle\!\!\!\Biggl\langle#2\Biggr\rangle\!\!\!\Biggr\rangle_{{#1}}}

\newcommand{\negquad}{\hskip-.5em\relax}

\begin{document}

\title{Feynman diagrams and the KdV hierarchy}
\author{Domenico Fiorenza}
\address{Dipartimento di Matematica ``Guido Castelnuovo'' ---
Universit\`a degli Studi di Roma ``la Sapienza'' --- P.le Aldo Moro, 2 --
00185 -- Roma, Italy} 

\email{fiorenza@mat.uniroma1.it}

\begin{abstract}The generating series of the intersection numbers of the stable
cohomology classes on moduli spaces of curves satisfies the string
equation and a KdV hierarchy. Kontsevich's original proof of this result uses a matrix model and the matrix Airy equation. Witten then recasted Kontsevich's results
in terms of Virasoro algebras, by means of an ingenious mixture of
Feynman diagrams techniques and integrations by parts with some
``rather formidable choice'' of the integrands. In this note we show how Witten's formidable choices can be bartered for standard Feynman diagram manipulations, as soon as one suitably enlarges the class of Feynman diagrams occurring in the proof.
\end{abstract}


\maketitle

\section*{Introduction}
A classical result of Kontsevich-Witten states that the
generating series of the intersection numbers of the
stable cohomology classes on moduli spaces of curves
satisfies the string equation and a KdV hierarchy
(see
\cite{ea;sketches} for a very comprehensive introduction to
the subject). Kontsevich's original proof
\cite{kontsevich;intersection-theory;1992} uses a matrix model, now
called the 't~Hooft-Kontsevich model, and the matrix Airy equation.
Witten \cite{witten;kontsevich-model} recasted Kontsevich's results
in terms of Virasoro algebras. His proof is an ingenious mixture of
Feynman diagrams techniques and integrations by parts with some
``rather formidable choice'' of the integrands.\footnote{
For instance, to prove the last equation in his paper, Witten
uses the familiar fact that \(\displaystyle{\int\ud
Y\frac{\del}{\del Y_{ij}}\left(M_{ij} \exp^{\frac{1}{2}\tr
\Theta Y^2-\frac{1}{6}\tr Y^3}\right)=0}\), where \(Y\) is an
\(N\times N\) anti-Hermitian matrix, \(\Theta\) is a negative
definite symmetric \(N\times N\) real matrix and \(M\) is any
polynomial in
\(Y\) and \(\Theta\). Then he makes the rather formidable
choice
\begin{align*}
M=\Theta^4&Y-\Theta^3Y\Theta+\frac{1}{2}\Theta^2Y\Theta^2+
\Theta^3Y^2-\Theta^2Y\Theta Y+\frac{3}{2}\Theta
Y\Theta^2Y-\frac{1}{4}\Theta Y\Theta Y \Theta-\Theta^2
Y^2\Theta\\
&+\Theta^2Y^3-\Theta Y\Theta Y^2-\frac{3}{2}\Theta Y^2\Theta
Y +\frac{9}{8}Y\Theta Y\Theta Y+\frac{1}{2}Y\Theta^2 Y^2 +
\frac{3}{2}\Theta Y^4-\frac{11}{8}Y\Theta Y^3\\
&-\frac{1}{16}Y^2\Theta
Y^2+\frac{1}{32}Y^5-\frac{N}{8}Y^2+\frac{5N}{2}\Theta
Y+\frac{3N}{2}\Theta^2-2\tr\Theta^2-\frac{1}{16}\tr
Y^2-\frac{3}{8}\tr \Theta Y\\
&+\frac{1}{4}Y\tr
Y-\frac{5}{2}Y\tr\Theta+\frac{1}{2}\Theta\tr Y\,.
\end{align*}}
The aim of this note is to show how the techniques developed
in \cite{fiorenza-murri;matrix-integrals}  allow to
straightforwardly rewrite Witten's proof entirely in terms of
Feynman diagrams, trading the formidable choices of the integrands for standard Feynman diagram manipulations.
The trick used to accomplish this is to use a wider class of Feymnan diagrams than the one considered by Witten in \cite{witten;kontsevich-model}. Namely, the class of Feynman diagrams with $n$-valent vertices labelled  by polynomials in $n$ variables, as in  \cite{fiorenza-murri;matrix-integrals}.
\bigskip\par
\noindent The paper is organized as follows
\medskip\par
\noindent 1\quad Intersection numbers on the moduli space of
curves\dotfill\pageref{sec:intersection-numbers}
\par
\noindent 2\quad The KdV hierarchy and Virasoro
operators\dotfill\pageref{sec:virasoro}
\par
\noindent 3\quad The idea of the proof\dotfill\pageref{sec:idea}
\par
\noindent 4\quad The 't~Hooft-Kontsevich matrix
model\dotfill\pageref{sec:tHooft}
\par
\noindent 5\quad Witten's formula for
derivatives\dotfill\pageref{sec:witten}
\par
\noindent 6\quad Proof of equation
(\ref{eq:i})\dotfill\pageref{sec:eq-i}
\par
\noindent 7\quad Proof of equation
(\ref{eq:ii})\dotfill\pageref{sec:eq-ii}
\par
\noindent Appendix: Formal differential
operators\dotfill\pageref{appendix}

\section{Intersection numbers on the moduli space of
curves}\label{sec:intersection-numbers}

For fixed integers \(g\geq 0\) and
\(n\geq 1\) with \(2-2g-n < 0\), let
\(\M_{g,n}\) be the moduli space of smooth complete curves of genus
\(g\) with \(n\) marked points and
\(\Mbar_{g,n}\) be its Deligne-Mumford compactification
\cite{deligne-mumford} (see \cite{arbarello-cornalba-griffiths} for an overview). Denote by
$\psi_i\in H^2(\Mbar_{g,n}, \setC)$ the Miller-Witten
cohomology classes \cite{miller,witten;2dgravity}, with
$i=1,\dots,n$, and by
$\langle \tau_{\nu_1} \cdots \tau_{\nu_n}
\rangle_{g,n}$ the intersection number \cite{witten;2dgravity}
\begin{equation*}
  \langle \tau_{\nu_1} \cdots \tau_{\nu_n} \rangle_{g,n}:=
\int_{\overline{\mathcal
      M}_{g,n}} \psi_1^{\nu_1} \cdots \psi_n^{\nu_n}.
\end{equation*}
 The generating series of intersection numbers (``free energy
 functional'' in physics literature) is the formal series in the variables
\(t_0,t_1,\dots\)
  \begin{equation*}
    F(t_*)=\sum_{g, n}F_{g, n}(t_*) := \sum_{g, n}
    \left(\frac{1}{n!} \sum_{\nu_1, \dots,\nu_n}\langle\tau_{\nu_1} \cdots
      \tau_{\nu_n}\rangle_{g, n}^{}t_{\nu_1}\cdots
      t_{\nu_n}\right)\,.
  \end{equation*}
 The exponential of the generating series is called the
\emph{partition function} of the intersection numbers and it is
denoted by the symbol $Z(t_*)$.

 It is well-known that the moduli
spaces of stable curve have an ideal orbi-cellularization whose
cells are indexed by isomorphism classes of closed connected ribbon
graphs with numbered holes, see, e.g.,
\cite{harer;cohomology-of-moduli,kontsevich;intersection-theory;1992, mondello-combinatorial, mondello-survey}.
As a consequence of this
fact, any integral on a moduli space of stable curves can be written as
a sum over isomorphism classes of numbered ribbon graphs. In
particular, Kontsevich finds the following
remarkable identity (Kontsevich's Main Identity):
\begin{multline}
    \sum_{\nu_1, \ldots, \nu_n} \langle\tau_{\nu_1} \cdots
    \tau_{\nu_n}\rangle_{g,n} \prod_{i=1}^n \frac{(2\nu_i - 1)!!}
    {\lambda_i^{2\nu_i + 1}}=
    \\
  \label{eq:KMI}
  = \sum_{(\Gamma,h)} \frac{1} {\card{\Aut
(\Gamma,h)}}\left(
      \frac{1} {2}\right)^{\card{\text{Vertices}(\Gamma)}}
\!\!\!\!\!\!\!\prod_{l\in
     \text{Edges}(\Gamma)}
\frac{2}{\lambda_{h(l^+)} +
\lambda_{h(l^-)}}
  \end{multline}
  where: the $\lambda_i$ are positive real variables; $(\Gamma,h)$ ranges
  over the set of isomorphism classes of closed connected numbered
  ribbon graphs of genus $g$ with $n$ holes; for any edge
$l$ of $\Gamma$, $l^+, l^-$ denote the
  (not necessarily distinct) holes $l$ belongs to.

\section{The KdV hierarchy and Virasoro operators}\label{sec:virasoro}

Let \(\setC[u,\del u/\del t_0,\dots]\) be the differential
algebra generated by the variable \(u\).
The \emph{Gel'fand -Dikii} polynomials are the
differential polynomials \(R_i(u)\)  
defined by the recursion
\begin{equation*}
\frac{\partial R_{n+1}(u)}{\partial
t_0}=\frac{1}{2n+1}\left(\frac{\partial u}{\partial t_0}
+2u\frac{\partial }{\partial
t_0}+\frac{1}{4}\frac{\partial^3}{\partial t_0^3}\right) R_n(u)\,,
\end{equation*}
with initial datum
\(
 R_0(u)=1
\)
and boundary condition
\( R_{n}(0)=0\) for \(n>0\).

The Korteweg-de Vries (KdV for short) hierarchy is the following
hierarchy of
differential equations for an element \(U\) of \(\setC[[t_*]]\):
\begin{equation*}
\frac{\partial U}{\partial t_{i}}=\frac{\partial}{\partial t_0}
 R_{i+1}(U)\,.
\end{equation*}
The first equation of the KdV hierarchy is
\begin{equation}\label{eq:KdV-one}
\frac{\del U}{\del t_1}=U\frac{\del U}{\del t_0}+\frac{1}{12}\frac{\del^3
U}{\del {t_0}^3}\,.
\end{equation}
If we set
\(\varphi(x,t):=U(\frac{1}{\sqrt{2}}x,3\sqrt{2}\,t,0,0,\dots)\),  equation
\prettyref{eq:KdV-one} becomes
\begin{equation*}
\varphi_t=6\varphi\varphi^{}_x+\varphi^{}_{xxx}\,,
\end{equation*}
which is the classical KdV equation (see \cite{ea;sketches} or
\cite{arnold}).

We are now ready to state the result we are going to prove
in this paper.
\begin{theorem*}[Kontsevich-Witten]
Let  \(F(t_*)\) be the generating series of intersection
numbers on moduli spaces of curves. Then
\begin{enumerate}
\item the series \(F\) satisfies the \emph{string equation}
\begin{equation*}
\frac{\partial F}{\partial
t_0}=\sum_{i=0}^{\infty}t_{i+1}\frac{\partial F}{\partial
t_i}+\frac{{t_0}^2}{2}\,;
\end{equation*}
\item the series \(U=\del^2F/\del {t_0}^2\) satisfies the
KdV hierarchy.
\end{enumerate}
\end{theorem*}

In \cite{witten;kontsevich-model} Witten showed that this theorem can
be recasted by means of Virasoro operators as follows. For any
\(\rho\in\setZ+1/2\), define a differential operator
\(\alpha_\rho\),
 by
\begin{equation*}
\alpha_\rho=\left\{\begin{matrix}
\displaystyle{\frac{(2\rho!!)}{\sqrt{2}}\frac{\partial}{\partial
t_{\rho+\frac{1}{2}}}} &\text{if}&\rho>0\\
\\
\displaystyle{\frac{1}{(-2\rho-2)!!\sqrt{2}}(
t_{-\rho-\frac{1}{2}}-\delta_{\rho+\frac{3}{2},0})}&\text{if}&\rho<0
\end{matrix}\right..
\end{equation*}
The operators \(\alpha_\rho\) satisfy the commutation relation
\begin{equation*}
[\alpha_{\rho_1},\alpha_{\rho_2}]=\rho_1\delta_{\rho_1+\rho_2,0}.
\end{equation*}
This implies that the formal differential operators
\begin{equation*}
L_n=\left\{\begin{matrix}
\displaystyle{\frac{1}{2}\sum_\rho \alpha_\rho\alpha_{n-\rho}},& n\neq0\\
\\
\displaystyle{\sum_{\rho>0}\alpha_{-\rho}\alpha_\rho+\frac{1}{16}},& n=0
\end{matrix}\right.
\end{equation*}
realize a Virasoro algebra with central charge \((c^3-c)/12\):
\begin{equation}\label{eq:Virasoro}
[L_m,L_n]=(m-n)L_{m+n}+\delta_{m+n,0}\frac{m^3-m}{12}
\end{equation}

Recasted by means of the $L_n$'s, the
Kontsevich-Witten theorem becomes
\begin{theorem*}
The partition function $Z(t_*)$ is a zero vector
for the Virasoro operators $L_n$, $n\geq -1$.
\end{theorem*}
By the commutation relations \prettyref{eq:Virasoro}, to prove the
Kontsevich-Witten's theorem one only needs to check \(L_nZ=0\) for
\(n=-1\) and \(n=2\). Written out explicitly, these two
equations are:
\begin{align*}
\tag{KW1}\label{eq:kw1}
&\frac{\partial Z}{\partial t_0}=
\sum_{i=0}^\infty  t_{i+1}
\frac{\partial Z}{\partial
t_{i}}+\frac{{t_0}^2}{2}Z\\
\tag{KW2}\label{eq:kw2}
&\frac{\partial Z}{\partial
t_3}
 =\frac{1}{7!!}\left(\sum_{i=0}^\infty
(2i+5)(2i+3)(2i+1)  t_i
\frac{\partial Z}{\partial t_{i+2}}+3\frac{\partial^2 Z}{\partial
t_0\partial
t_1}\right)\,.
\end{align*}

\section{The idea of the proof}\label{sec:idea}
As a consequence of its Main Identity among the intersection numbers
\(\langle\tau_{\nu_1}\!\!\cdots\tau_{\nu_n}\!\rangle\), Kontsevich
proves that the partition function $Z(t_*)$ is related to the asymptotic
expansion of the Hermitian matrix integral
\begin{equation}\label{eq:hermitean-integral}
\int_{\Hermitian[N]}\exp\left\{\frac{\U}{6}\tr X^3\right\}
\ud\mu_\Lambda(X)\,,
\end{equation}
where
$\Lambda=\diag\{\Lambda_1,\dots,\Lambda_N\}$ is
a positive definite Hermitian matrix and \(\ud\mu_\Lambda(X)\) is the
Gaussian measure defined by the inner product
\begin{equation*}
{\inner{X}{Y}}_\Lambda:=\onehalf \left(\tr(X\Lambda Y) +
\tr(Y\Lambda X)\right)
\end{equation*}
 on the space
\(\Hermitian[N]\) of $N\times N$ Hermitian matrices. More precisely,
let, for any $k\geq 0$,
\begin{equation*}
t_k(\Lambda)=-(2k-1)!!\tr\Lambda^{-(2k+1)}\,.
\end{equation*}
Then
\begin{equation}
Z(t_*)\biggr\vert_{t_*(\Lambda)}\asymp
\int_{\Hermitian[N]}\exp\left\{\frac{\U}{6}\tr X^3\right\}
\ud\mu_\Lambda(X)
\end{equation}
as $|\Lambda|\to\infty$. When \(N\to\infty\), the
\(\{t_k(\Lambda)\}\) become independent coordinates (Miwa
coordinates) on the space
\(\Hermitian[N]/U(N)\), so, the
Kontsevich-Witten's theorem is equivalent to
\begin{equation*}
\left.L_nZ\right\vert_{t_*(\Lambda)}=0,\quad n=-1,2,
\end{equation*}
i.e., it is reduced to a statement concerning the Hermitian
matrix integral \prettyref{eq:hermitean-integral}. One of the big
advantages of this translation of the original problem in a problem
concerning Hermitian matrices is that the differential operators
\begin{equation*}
\sum_{i=0}^\infty  t_{i+1}
\frac{\partial}{\partial
t_{i}}\qquad \text{and}\qquad
\sum_{i=0}^\infty
(2i+5)(2i+3)(2i+1)  t_i
\frac{\partial}{\partial t_{i+2}}
\end{equation*}
appearing in equations (\ref{eq:kw1}-\ref{eq:kw2}) are very natural
in terms of the Miwa coordinates. Indeed they
correspond to the action of the operators
$\tr\Lambda^{-1}\partial/\partial\Lambda$ and
$\tr\Lambda^{5}\partial/\partial\Lambda$ on a function of the
$t_*(\Lambda)$. More precisely, we have\footnote{
Both equations immediately follow by
\begin{equation*}
\tr\Lambda^{2k+1}\frac{\partial}{\partial
\Lambda}t_i(\Lambda)=
-\frac{(2i+1)!!}{(2i-2k-1)!!}t_{i-k}(\Lambda)\,,
\end{equation*}
for any \(k\geq-1\) and any \(i\geq
k\).
}
\begin{align*}
\tag{KW1'}\label{eq:kw1'}
&\tr\Lambda^{-1}\frac{\partial}{\partial
\Lambda}Z(t_*(\Lambda))=-\left(\sum_{i=0}^{\infty}\left.
t_{i+1}\frac{\partial
Z}{\partial t_i} \right)\right   
\rvert_{t_*(\Lambda)}\\
&\tr\Lambda^{5}\frac{\partial}{\partial
\Lambda}Z(t_*(\Lambda))=\left(-\sum_{i=0}^{\infty}
(2i+5)(2i+3)(2i+1)
t_{i}\frac{\partial Z}{\partial t_{i+2}}\right.+\\
\tag{KW2'}\label{eq:kw2'}&\phantom{mmmmmmmmmmmmmmm}
\left.\left.
+3\tr\Lambda\frac{\del Z}{\del
t_1}+\tr
\Lambda^{3}\frac{\del Z}{\del t_0}
\right)\right   
\rvert_{t_*(\Lambda)}
\end{align*}

Put together equations
(\ref{eq:kw1}-\ref{eq:kw2}) and equations (\ref{eq:kw1'}-\ref{eq:kw2'})
to obtain that the Kontsevich-Witten's theorem is equivalent to the
following two equations:
\begin{align*}
\tag{I}\label{eq:i}
&\left.\frac{\partial Z}{\partial t_0}\right\vert_{t_*(\Lambda)}=
-\tr\Lambda^{-1}\frac{\del }{\del
\Lambda}Z(t_*(\Lambda))+\frac{(\tr
\Lambda^{-1})^2}{2}Z(t_*(\Lambda))\\
\tag{II}\label{eq:ii}
&\left.\left(-105\frac{\partial Z}{\partial
t_3}+3\frac{\partial^2 Z}{\partial{t_0}\partial
t_1}+3\tr\Lambda\frac{\del Z}{\del
t_1}+\tr\Lambda^3\frac{\del Z}{\del
t_0}\right)\right\vert_{t_*(\Lambda)}=
\tr\Lambda^{5}\frac{\del }{\del
\Lambda}Z(t_*(\Lambda))
\end{align*}

 Having written
\(Z(t_*(\Lambda))\) as an Hermitian matrix integral it is immediate
to compute that also \((\tr\Lambda^{-1}\del/\del
\Lambda)Z(t_*(\Lambda))\) and \((\tr\Lambda^{5}\del/\del
\Lambda)Z(t_*(\Lambda))\) are Hermitian matrix integrals. Moreover,
one can show that
 also the derivatives of the partition function
$Z(t_*)$ with respect to the
\(t_*\) variables (evaluated at $t_*(\Lambda)$) can be expressed
as Hermitian matrix integrals.\footnote{ In
\cite{witten;kontsevich-model} Witten shows that for any
\(D\) in \(\{\del/\del t_0,\,
\del/\del t_1,\, \del/\del t_2,\, \del/\del t_3,\,\del^2/\del {t_0}^2,\break
\del^2/\del {t_0}\del t_1\}\), there exist a polynomial
\(P_D\)  in the odd traces of \(X\) such that
\begin{equation*}
DZ(t_*)\biggr\vert_{t_*(\Lambda)}=\int_{\Hermitian[N]}
P_D(X)\cdot\exp\left\{\frac{\U}{6}\tr X^3\right\}
\ud\mu_\Lambda(X)\,,
\end{equation*}
and conjectures that this should  be true for any
differential operator in the variables \(t_*\). The Witten
conjecture has been proved by Di~Francesco, Itzykson and
Zuber in \cite{di-francesco-itzykson-zuber;kontsevich-model}
and is henceforth known as the DFIZ theorem.
We proved a generalization of the DFIZ theorem in
\cite{fiorenza-murri;matrix-integrals}. See \cite{arbarello-cornalba;dfiz, bini} for the use of the 
DFIZ theorem in the investigation of the geometry of moduli spaces of pointed curves.}
 The proof of
equations (\ref{eq:i}-\ref{eq:ii}) is therefore reduced to
checking two identities between Hermitian matrix integrals.

Such a check is performed by Witten \cite{witten;kontsevich-model} by a
clever mixture of Feynman diagrams techniques and tricky integration by
parts. We are going to show how introducing suitable Feynman
diagrams with vertices decorated by polynomials as in
\cite{fiorenza-murri;matrix-integrals} reduces everything to
checking a few graphical identities.

The rest of this paper is organized as follows: we first
recall some general notation and result concerning the
relation between Feynman diagrams and Gaussian integrals.
Next we introduce the Feynman rules for the
't~Hooft-Kontsevich model and rewrite the right-hand sides
of equations (\ref{eq:i}-\ref{eq:ii}) as expectation values of
suitable Feynman diagrams. We will then discuss Witten's
formula for the derivatives of the partition function
$Z(t_*)$, which allows to write also the left-hand sides of
equations (\ref{eq:i}-\ref{eq:ii}) as expectation values of
suitable Feynman diagrams (with vertices decorated by
polynomials).  Having done this the proof of the
theorem is reduced to checking two identities
among expectation values of Feynman diagrams:
this is a straightforward (but quite long and tedious) computation.

\section{The 't~Hooft-Kontsevich matrix model}\label{sec:tHooft}

Here we will briefly recall the Feynman rules\footnote{
Details on these Feynman rules can be found in
\cite{fiorenza-murri;matrix-integrals}; we refer the reader to
\cite{fiorenza-murri;feynman-diagrams, fiorenza;integration-groupoids} and references therein for an
introduction to the general theory of Feynman diagrams and their
relations with the asymptotic expansion of Gaussian integrals.} for
the Gaussian integral
\begin{equation}\label{eq:gaussian-integral}
\int_{\Hermitian[N]}\exp\biggl\{\frac{\U}{6}\tr X^3\biggr\}
\ud\mu_\Lambda(X)
\end{equation}
where $\ud\mu_\Lambda(X)$ is the Gaussian measure on
$\Hermitian[N]$ induced by the inner product $(X|Y)_\Lambda=
\frac{1}{2}(\tr X\Lambda Y+\tr Y\Lambda X)$, i.e.,
\begin{equation*}
\ud\mu_\Lambda(X):=\frac{\ud X}{\displaystyle{
\int_{\Hermitian[N]}\exp\biggl\{-\frac{1}{2}\tr
\Lambda X^2\biggr\}
\ud X}}\,.
\end{equation*}

The space of fields is the complexification of the space
$\Hermitian[N]$ of $N\times N$ Hermitian matrices, so it is
naturally isomorphic to the space
$M_N(\setC)$ of
$N\times N$ complex matrices. We will denote by $\{E_{ij}\}$
the canonical basis of $M_N(\setC)$ and will write the
indices $ij$ near to a leg in a Feynman diagram to denote
evaluation or coevaluation at the basis element $E_{ij}$.

The graphical element corresponding to the field $M\in
M_N(\setC)$ will be denoted by
\begin{equation*}
{\xy,(0,-.5)*+[F]{M};(0,1)**\dir{-}
,(-.2,1)*{i},(.2,1)*{j}
\endxy}\,\mapsto M_{ij}
\end{equation*}
In particular, we have
\begin{equation*}
{\xy,(0,-.5)*+[F]{\Lambda^l};(0,1)**\dir{-}
,(-.2,1)*{i},(.2,1)*{i}
\endxy}\,\mapsto \Lambda_{i}^l
\end{equation*}

The \emph{propagator}
corresponds to the copairing dual to the inner product
$(-|-)_\Lambda$, so it is given by
\begin{equation*}
  {\xyc\vloop-
    ,(0.05,0.05)*\txt{\({}_i\ {}_j\)},(1,0.05)*\txt{\({}_j\ {}_i\)}
    \endxyc}
  \mapsto \frac {2} {\Lambda_i + \Lambda_j}
\end{equation*}
The \emph{interaction} corresponds to the cyclically
invariant tensor
\begin{equation*}
X\otimes Y\otimes Z\mapsto \frac{\U}{2}\tr
XYZ\,,
\end{equation*}
 so it is represented by a trivalent vertex
\begin{equation*}
  {\xyc%
    \rgvertex3%
    \fence1{i}{j}%
    \fence2{j}{k}%
    \fence3{k}{i}%
    \endxyc}
\mapsto \frac{\U}{2}
\end{equation*}
We well call this an \emph{ordinary} vertex, in order to
distinguish it from the \emph{special} vertices that we will
introduce below. Since the interaction enjoys a cyclic
invariance, the trivalent vertex representing it is equipped
with a cyclic order on the legs; therefore the Feynman
diagrams related to the Gaussian integral
\prettyref{eq:gaussian-integral} are
\emph{ribbon graphs}. In the displayed diagrams of this
paper, the cyclic order on the vertices will be the one
induced by the standard orientation of the plane.

We conclude the list of the Feynman rules of the
't~Hooft-Kontsevich model by introducing the special
vertices; they correspond
to the cyclic interactions
\begin{equation*}
X_1\otimes\cdots X_n\mapsto \tr(X_1\cdots X_n)
\end{equation*}
and are graphically represented as
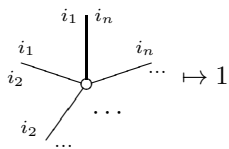
\begin{equation*}
  \xyc
  \rgvertex[\!]{5}\genericfences
  \endxyc
  \mapsto 1
\end{equation*}

Closed trivalent ribbon graph corresponding to
top-dimensional cells in the ideal orbi-cellularization of
the moduli spaces of curves can be considered as Feynman
diagrams in the 't~Hooft-Kontsevich model with only ordinary
vertices and no legs. By the Feynman rules described above,
is immediate to compute that the amplitude of such a ribbon
graph
$\Gamma$ is given by
\begin{equation}\label{eq:like-KMI-rhs}
 Z_{\Lambda}(\Gamma)=(-1)^{|\Holes{\Gamma}|}\left(\frac{1}{2}
\right)^{|\Vertices{\Gamma}|}\sum_{c}\prod_
{l\in\Edges{\Gamma}}\frac{2}{\Lambda_{c(l^+)}
+\Lambda_{c(l^-)}},
\end{equation}
where $\Holes{\Gamma}$ denotes the set of the holes of the
ribbon graph
$\Gamma$, $\Vertices{\Gamma}$ is set of the vertices of
$\Gamma$, $c$ ranges in the set of all maps
$\Holes{\Gamma}\to\{1,\dots,N\}$, and $l^\pm$ are the two
(not necessarily distinct) holes $l$ belongs to. So the
graphical elements we are considering are actually ribbon
graphs whose holes are coloured by the colours
$\{1,\dots,N\}$.

The right-hand side of equation \prettyref{eq:like-KMI-rhs} is very
similar to the right-hand side of the Kontsevich' Main
Identity \prettyref{eq:KMI}. Indeed, a bit of combinatorics
shows that
\begin{equation}\label{eq:cite-it-just-below}
Z(t_*)\biggr\rvert_{t_*(\Lambda)}=\sum_{\Gamma}\frac{Z_\Lambda(\Gamma)}{
\card{\Aut\Gamma}}
\end{equation}
where $\Gamma$ ranges into the set of closed trivalent ribbon
graphs. See \cite{kontsevich;intersection-theory;1992} or
\cite{fiorenza-murri;matrix-integrals} for details. The right
hand-side of equation \prettyref{eq:cite-it-just-below} is
the expectation value of the vacuum in the
't~Hooft-Kontsevich model. To see this, recall that
the \emph{expectation value} of the ribbon
graph $\Psi$ is defined as the sum
  \begin{equation*}
    \langle\!\langle \Psi\rangle\!\rangle_{\Lambda}^{}
    := \sum_{\Gamma \in \R_{\Psi}^{}(0)}
    \frac{Z_{\Lambda}^{}(\Gamma)}{\card{\Aut{\Gamma}}}\,.
  \end{equation*}
where  \(\R_{\Psi}^{}(0)\) denotes the set of (isomorphism
classes of) closed trivalent ribbon graphs containing
\(\Psi\) as a distinguished sub-graph and having no special
vertex outside $\Psi$. By saying that the sub-graph
\(\Psi\) is distinguished, we require that any automorphism
of an object
\(\Gamma \in\R_{\Psi}^{}(0)\) maps \(\Psi\) onto itself.  It
follows from the definition that \(\R_\emptyset^{}(0)\) is
the set of isomorphism classes of all closed trivalent ribbon
graphs.

If
$\Psi$ is a ribbon graph (possibly with special vertices)
with
$n$ legs, then
\begin{equation*}
X\mapsto Z_\Lambda(\Psi)(X^{\otimes n})
\end{equation*}
is a polynomial function (of degree $n$) on $\Hermitian[N]$
so it is integrable with respect to the Gaussian measure
$\ud\mu_\Lambda$ and,
by the general theory of Feynman diagrams, we have:
\begin{equation*}
   \langle\!\langle \Psi\rangle\!\rangle_{\Lambda}^{} =
\int_{\Hermitian[N]}
\frac{Z_{\Lambda}
      (\Psi)(X\tp{n})}{\card{\Aut\Psi}}
    \exp\biggl\{\frac{\U}{6}\tr X^3\biggr\}
\ud\mu_\Lambda(X)
  \end{equation*}

In particular, we have:
\begin{equation*}
  \langle\!\langle \emptyset\rangle\!\rangle_{\Lambda}^{}=
  \int_{\Hermitian[N]}\exp\biggl\{\frac{\U}{6}\tr X^3\biggr\}
\ud\mu_\Lambda(X)
\end{equation*}
so that equation \prettyref{eq:cite-it-just-below} can be rewritten as
\begin{equation*}
 Z(t_*)\biggr\vert_{t_*(\Lambda)} =
\int_{\Hermitian[N]}\exp\biggl\{\frac{\U}{6}\tr X^3\biggr\}
\ud\mu_\Lambda(X)\,.
\end{equation*}

We also introduce the subset $\overline{\Gamma}$ of
$\R_\Gamma(0)$ consisting of those Feynman diagrams having
no vertices except those in $\Gamma$, i.e., the set of
Feynman diagrams that can be obtained by joining the legs of
$\Gamma$ by means of edges. We call elements in
$\overline{\Gamma}$ the closures of $\Gamma$. Clearly, if
$\Gamma$ has an odd number of legs, then $\overline{\Gamma}$
is empty. By definition, the \emph{expectation value without
interactions} of the ribbon graph $\Psi$ is the
sum
  \begin{equation*}
\langle\!\langle \Psi\rangle\!\rangle_{\Lambda,0}^{}
    := \sum_{\Gamma \in \overline{\Psi}}
    \frac{Z_{\Lambda}^{}(\Gamma)}{\card{\Aut{\Gamma}}}\,.
  \end{equation*}
and we have
\begin{equation*}
    \langle\!\langle
\Psi\rangle\!\rangle_{\Lambda,0}^{} =
\int_{\Hermitian[N]}
\frac{Z_{\Lambda}
      (\Psi)(X\tp{n})}{\card{\Aut\Psi}}
    \ud\mu_\Lambda(X)
  \end{equation*}
Note that, if $\Psi$ is a closed ribbon graph, then
$\overline{\Psi}=\{\Psi\}$ so
\begin{equation*}
   \langle\!\langle \Psi\rangle\!\rangle_{\Lambda,0}^{}
    =
    \frac{Z_{\Lambda}^{}(\Psi)}{\card{\Aut{\Psi}}}\,,
\end{equation*}
and
\begin{equation*}
    \langle\!\langle \Psi\rangle\!\rangle_{\Lambda}^{}
    =
    \langle\!\langle \Psi\rangle\!\rangle_{\Lambda,0}^{}\cdot
\langle\!\langle \emptyset\rangle\!\rangle_{\Lambda}^{}\,.
\end{equation*}

We now come back to the two equations in the Kontsevich-Witten theorem.
It is easy to compute, for any \(l\in\setZ\),
\begin{align}
\notag
&\left(\tr\Lambda^l\frac{\partial}{\partial
\Lambda}\right)\int_{\Hermitian[N]}\exp\left\{\frac{\U}{6}\tr
X^3\right\}\ud\mu_\Lambda(X)=\\
\notag
&\qquad=
-\frac{1}{2}\int_{\Hermitian[N]}\tr \Lambda^l
X^2\exp\left\{\frac{\U}{6}\tr
x^3\right\}\ud\mu_\Lambda(X)+\\
\label{eq:detrlambda}&
\qquad\qquad+\frac{1}{2}\int_{\Hermitian[N]}\exp\left\{\frac{\U}{6}\tr
X^3\right\}\ud\mu_\Lambda(X)\cdot
\int_{\Hermitian[N]}\tr \Lambda^l X^2\ud\mu_\Lambda(X)
\end{align}
The bilinear map
\begin{equation*}
X\otimes Y\mapsto \tr X\Lambda^l Y
\end{equation*}
can be seen as the amplitude of the diagram
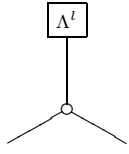
\begin{equation*}
\xy
\rgvertex[\!]{3}\loose1\loose2\loose3%
,(0,1.7)*+[F]{\Lambda^l};(0,.5)**\dir{-}
\endxy\end{equation*}
so that equation \prettyref{eq:detrlambda} can be rewritten
as
\begin{align*}
\left(\tr\Lambda^l\frac{\partial}{\partial
\Lambda}\right)Z(t_*(\Lambda))&=-\frac{1}{2}\gint{
{\xy\rgvertex[{\!}]{3}\loose1\loose2\loose3%
,(0,1.7)*+[F]{\Lambda^l};(0,.5)**\dir{-}
\endxy}
}
+\frac{1}{2}\langle\!\langle\emptyset\rangle\!\rangle_\Lambda^{}
\cdot\Biggl\langle
{\xy\rgvertex[{\!}]{3}\loose1\loose2\loose3%
,(0,1.7)*+[F]{\Lambda^l};(0,0.5)**\dir{-}
\endxy}
\Biggr\rangle_0\\
\notag&=-\frac{1}{2}\gint{
{\xy\rgvertex[{\!}]{3}\loose1\loose2\loose3%
,(0,1.7)*+[F]{\Lambda^l};(0,.5)**\dir{-}
\endxy}
}
+\frac{1}{2}\langle\!\langle\emptyset\rangle\!\rangle_\Lambda^{}
\cdot Z_\Lambda\Biggl(
{\xy\rgvertex[{\!}]{3}\loose1\join23%
,(0,1.7)*+[F]{\Lambda^l};(0,.5)**\dir{-}
\endxy}
\Biggr)
\end{align*}

It immediate to compute
\begin{align*}
Z_\Lambda\left({\xy\rgvertex[\!]{3}\loose1\join23%
,(0,1.7)*+[F]{\Lambda^l};(0,.5)**\dir{-}%
\endxy}\right)&=\sum_{i,j}
{\xy\rgvertex[\!]{3}\loose1\join23%
,(0,1.7)*+[F]{\Lambda^l};(0,.5)**\dir{-}%
,(0,-.3)*{j},(-.2,0.3)*{i}
\endxy}=\sum_{i,j}\frac{2{\Lambda_i}^l}{(\Lambda_i+\Lambda_j)}=
\sum_{i,j}\frac{{\Lambda_i}^l+
{\Lambda_j}^l}{
\Lambda_i+\Lambda_j}
\end{align*}
So, in particular
\begin{equation*}
Z_\Lambda\left({\xy\rgvertex[\!]{3}\loose1\join23%
,(0,1.7)*+[F]{\Lambda^{-1}};(0,.5)**\dir{-}%
\endxy}\right)=
(\tr\Lambda^{-1})^2
\end{equation*}
and
\begin{equation*}
Z_\Lambda\left({\xy\rgvertex[\!]{3}\loose1\join23%
,(0,1.7)*+[F]{\Lambda^{5}};(0,.5)**\dir{-}%
\endxy}\right)=
2\tr\Lambda^{4}-2\tr\Lambda^3\tr\Lambda+(\tr\Lambda^2)^2
\end{equation*}
Therefore, equations (\ref{eq:i}-\ref{eq:ii}) can be rewritten as
\begin{align*}
\tag{I}
&\left.2\frac{\partial Z}{\partial t_0}\right\vert_{t_*(\Lambda)}=
\gint{
{\xy\rgvertex[{\!}]{3}\loose1\loose2\loose3%
,(0,1.7)*+[F]{\Lambda^{-1}};(0,.5)**\dir{-}
\endxy}
}
\\ \\
\tag{II}
&\left.\left(210\frac{\partial Z}{\partial t_3}-6\frac{\partial^2
Z}{\partial{t_0}\partial
t_1}+6\tr\Lambda\frac{\del Z}{\del
t_1}+2\tr\Lambda^3\frac{\del Z}{\del
t_0}\right)\right\vert_{t_*(\Lambda)}=
\\
&\phantom{mmm}=
\gint{
{\xy\rgvertex[{\!}]{3}\loose1\loose2\loose3%
,(0,1.7)*+[F]{\Lambda^5};(0,.5)**\dir{-}
\endxy}
}-\left(2\tr\Lambda^{4}-2\tr\Lambda^3\tr\Lambda+(\tr\Lambda^2)^2
\right)\cdot\langle\!\langle\emptyset\rangle\!\rangle_\Lambda^{}
\end{align*}

\section{Witten's formula for derivatives}\label{sec:witten}

We have seen in \prettyref{sec:tHooft} that the Kontsevich'
Main Identity relates the partition function of the
intersection numbers $Z(t_*)$ to the expectation value of
the vacuum in the $N$-dimensional 't~Hooft-Kontsevich model.
It has been remarked by Witten
\cite{witten;kontsevich-model} that the Kontsevich'
Main Identity also relates the first order derivatives
of the partition function $Z(t_*)$ to
expectation values in the $(N+1)$-dimensional
't~Hooft-Kontsevich model of ribbon graph with one
``distinguished'' hole, and more generally, the
$k^{\text{th}}$ order derivatives of $Z(t_*)$ to to
expectation values in the $(N+k)$-dimensional
't~Hooft-Kontsevich model of ribbon graph with $k$
``distinguished'' holes.

Indeed, by the definition of the free energy $F(t_*)$ it
immediately follows that
\begin{equation*}
  \frac {\del F(t_*)} {\del t_k} = \sum_{g,\nu_*}
  \frac{1}{n!}\langle
  \tau_{\nu_1} \ldots \tau_{\nu_n}\tau_k  \rangle_{g,n}
  t_{\nu_1} \cdots t_{\nu_n},
\end{equation*}
so that taking the derivative with respect to the variable
$t_k$ corresponds to ``adding a $\tau_k$ at the end of
the intersection indices''. If we write the Kontsevich' Main
Identity for ribbon graphs with $(n+1)$ holes, we find
 \begin{multline*}
   \sum_{k} \left(\sum_{\nu_1, \dots, \nu_n} \langle
      \tau_{\nu_1} \cdots \tau_{\nu_n}
\tau_{k}\rangle_{g,n} \frac{{(2 \nu_i-1)!!}}{
{\lambda_{i}^{2 \nu_i+1}}}
   \right) \frac{{(2 k-1)!!}}{ {\lambda_{n+1}^{2 k+1}}}=
    \\
    =\sum_{ (\Gamma,h)} \frac{1}{\card{\Aut{(\Gamma,h)}}}
    \left( \frac{1}{2}\right)^{|\Vertices{\Gamma}|}
    \prod_{l\in\Edges{\Gamma}} \frac{2}{\lambda_{h(l^+)} +
      \lambda_{h(l^-)}},
  \end{multline*}
where $(\Gamma,h)$ ranges in the set of isomorphism classes
 of closed
connected trivalent genus $g$ ribbon graphs,
 with
$n+1$ holes, numbered from $1$ to $n+1$.  Therefore,
\begin{multline*}
  \sum_{m_*,\nu_*} \langle
  \tau_{\nu_1}\cdots\tau_{\nu_n} \tau_{k}\rangle_{g,n+1}
  \prod_{i=1}^n\frac{(2 \nu_i-1)!!}{\lambda_{i}^{2n_i+1}} =
  \\
  =\frac{1}{(2k-1)!!} \Coeff_{\lambda_{n+1}}^{-(2k + 1)} \Biggl( \sum_{
(\Gamma,h)}
  \frac{(1/2)^{|\Vertices{\Gamma}|}}{\card{\Aut{(\Gamma,h)}}}
\prod_{l\in\Edges{\Gamma}}
  \frac{2}{\lambda_{h(l^+)}+ \lambda_{h(l^-)}} \Biggr).
\end{multline*}

Therefore, the derivatives of the generating series $F(t_*)$
can be seen as Laurent coefficients of the Main Identity
with one more variable $\lambda_{n+1}^{}$. Translating this
fact into the language of Hermitian matrix integrals lead us
to consider the $(N+1)$-dimensional 't~Hooft-Kontsevich
model
$(N+1)$-dimensional 't~Hooft-Kontsevich model
 $Z_{z\oplus\Lambda}$ defined by the
diagonal matrix
\begin{equation*}
  z\oplus \Lambda =
  \begin{pmatrix}
    z       & 0\\
    0 & \Lambda
  \end{pmatrix},
\end{equation*}
where
\(z\) is a positive real variable.
The propagators in
the \((N+1)\)-dimensional `t~Hooft-Kontsevich model are
\begin{align}
  \label{eq:casimirs-n-piu-uno1}
  {\xyc\vloop-
    ,(0.05,0.05)*\txt{\({}_i\ {}_j\)},(1,0.05)*\txt{\({}_j\ {}_i\)}
    \endxyc}
  \mapsto \frac {2} {\Lambda_i + \Lambda_j}
  & &
  {\xyc\vloop-
    ,(0.05,0.05)*\txt{\({}_z\ {}_i\)},(1,0.05)*\txt{\({}_i\ {}_z\)}
    \endxyc}
  \mapsto \frac {2} {z + \Lambda_i} &
\\
\notag
\\
\label{eq:casimirs-n-piu-uno2}
 {\xyc\vloop-
    ,(0.05,0.05)*\txt{\({}_i\ {}_z\)},(1,0.05)*\txt{\({}_z\ {}_i\)}
    \endxyc}
  \mapsto \frac {2} {z+ \Lambda_i}
& &
  {\xyc\vloop-
    ,(0.05,0.05)*\txt{\({}_z\ {}_z\)},(1,0.05)*\txt{\({}_z\ {}_z\)}
    \endxyc}
  \mapsto \frac {1} {z} &
\end{align}
so the graphical elements we are considering are actually
ribbon graphs with holes coloured by the colours
$\{z,1,2,\dots,N\}$. We will put a \(z\) in the middle of an
hole to mean that the hole is given the ``special'' colour
$z$; colours on a blank hole will be allowed to vary in the
set $\{1,2,\dots,N\}$. A little combinatorics (see details
in \cite{fiorenza-murri;matrix-integrals}) then shows that,
for any $k\geq0$, the following identity holds:
\begin{equation}\label{eq:dZst}
  \left.\frac {\del Z(t_*)} {\del t_k}\right\rvert_{t_*(\Lambda)}
  = -\frac{1}{(2k-1)!!}  
  \Coeff_{z}^{-(2k + 1)}\! \Biggl(
    \sum_{\Gamma\in\R_\emptyset^{[1]}(0)}
    \frac{Z_{z\oplus\Lambda}(\Gamma)} {\card{\Aut{\Gamma}}}
    \Biggr),  
\end{equation}
where $\R_\emptyset^{[1]}(0)$ denotes the set of isomorphism classes
of trivalent closed ribbon graphs with  exactly
one $z$-decorated hole (and all the other holes blank).

A completely similar argument works for higher order derivatives. For
instance,
second order derivatives are related to amplitudes in the
$(N+2)$-dimensional
't~Hooft-Kontsevich model defined by the diagonal matrix
\begin{equation*}
  w\oplus z\oplus \Lambda =
  \begin{pmatrix}
    w &0 &0\\
    0 &z       & 0\\
    0 & 0 & \Lambda
  \end{pmatrix}.
\end{equation*}
by the formula
\begin{equation}\label{eq:dZst-higher}
  \left.\frac {\del^2 Z(t_*)} {\del t_{k_1}\del
t_{k_2}}\right\rvert_{t_*(\Lambda)}
  = \frac{\Coeff_{w}^{-(2k_1 + 1)}\Coeff_{z}^{-(2k_2 +
1)}}{(2k_1-1)!!(2k_2-1)!!}  
  \Biggl(
    \sum_{\Gamma\in\R_\emptyset^{[1,1]}(0)}
    \frac{Z_{z\oplus\Lambda}(\Gamma)} {\card{\Aut{\Gamma}}}
    \Biggr),  
\end{equation}
where $\R_\emptyset^{[1,1]}(0)$ denotes the set of isomorphism classes
of trivalent closed ribbon graphs with only ordinary
vertices and exactly one $z$-decorated hole and one
$w$-decorated one (and all the other holes blank).

Note that we are extracting the Laurent coefficient of $z^{-(2k_2+1)}$
first, and
the coefficient of $w^{-(2k_1+1)}$ later; this means that we are
considering the
Laurent expansion in the region $|z|>|w|$. It follows from the
Kontsevich' Main
Identity that one would get the same result by considering the
Laurent expansion in the region $|w|>|z|$.

Now, let us give a closer look to equation
\prettyref{eq:dZst}. If $\Gamma$ is an element of
$\R_\emptyset^{[1]}$, then its
$z$-decorated hole can be regarded as a distinguished
sub-diagram of $\Gamma$. So, if we call ($z$-decorated)
hole-type a minimal element of $\R_\emptyset^{[1]}$, i.e.,
a ribbon graph $\Gamma\in\R_\emptyset^{[1]}$ such that no
proper subgraph of $\Gamma$ lies in $\R_\emptyset^{[1]}$,
then equation \prettyref{eq:dZst}
can be rewritten as:
\begin{equation*}
  \left.\frac {\del Z(t_*)} {\del t_k}\right\rvert_{t_*(\Lambda)}
  = -\frac{1}{(2k-1)!!}\sum_{\Gamma\in{\mathcal S}_z}
  \Coeff_{z}^{-(2k + 1)}\! \Biggl(
    \sum_{\Phi\in\R_\Gamma^{[1]}(0)}
    \frac{Z_{z\oplus\Lambda,s_*}(\Phi)} {\card{\Aut{\Phi}}}
    \Biggr)\,,
\end{equation*}
where ${\mathcal S}_z$ denotes the set of
isomorphism classes of hole types, and $\R_\Gamma^{[1]}(0)$
denotes the set of isomorphism  classes of closed ribbon
graphs containing the hole type
$\Gamma$ as a distinguished subgraph and having no
$z$-decorated hole apart from the hole of $\Gamma$.
 By
introducing the shorthand notation
\begin{equation*}
  \langle\!\langle{\Gamma}\rangle\!\rangle_{z\oplus\Lambda}^{[1]}
  := \sum_{\Phi\in\R^{[1]}_\Gamma(0)}
  \frac{Z_{z\oplus\Lambda}^{}(\Phi)} {\card{\Aut\Phi}}\,,
\end{equation*}
where $\Gamma$ is an hole type, Witten's formula for derivatives
finally becomes:
\begin{equation}\label{eq:dZst-twotwo}
  \left.\frac {\del Z(s_*;t_*)} {\del t_k}\right\rvert_{t_*(\Lambda)}
  = -\frac{1}{(2k-1)!!}\sum_{\Gamma\in{\mathcal S}_z}
  \Coeff_{z}^{-(2k +
1)}\langle\!\langle{\Gamma}\rangle\!\rangle_{z\oplus\Lambda}^{[1]}
 \end{equation}

The sum on the right-hand side of equation
\prettyref{eq:dZst-twotwo} is actually a \emph{finite} sum.
Indeed, according to the Feynman rules
(\ref{eq:casimirs-n-piu-uno1}-\ref{eq:casimirs-n-piu-uno2}), each edge
of the
$z$-decorated hole brings a factor of order $O(z^{-1})$ as $z\to\infty$.
Then, the coefficient of $z^{-(2k+1)}$ in the Laurent
expansion of the amplitude $Z_{z\oplus \Lambda}(\Gamma)$
will be zero as soon as the hole-type $\Gamma$ has more than
$2k+1$ edges.
 Since there are finitely many hole types such that the
$z$-decorated hole is bounded by at most $2k+1$ edges, we
are actually dealing with a finite sum.

 Similarly, we call
($\{z,w\}$-decorated) two-holes-types the minimal elements
of $\R_\emptyset^{[1,1]}$, so that we can rewrite
equation \prettyref{eq:dZst-higher}
as
\begin{equation}\label{eq:dZst-twotwoo}
  \left.\frac {\del^2 Z(s_*;t_*)} {\del t_{k_1}\del
t_{k_2}}\right\rvert_{t_*(\Lambda)}
  = \sum_{\Gamma\in{\mathcal S}_{z,w}}\frac{\Coeff_{w}^{-(2k_1 +
1)}\Coeff_{z}^{-(2k_2 + 1)}}{(2k_1-1)!!(2k_2-1)!!}
  \langle\!\langle{\Gamma}\rangle\!\rangle_{w\oplus z\oplus\Lambda}^{[1,1]}
 \end{equation}
where the symbol ${\mathcal S}_{z,w}$ denotes
the set of isomorphism classes of two-holes types.
An argument completely similar to the one used above shows
that the also sum on the right-hand side of
\prettyref{eq:dZst-twotwoo} is actually a finite sum. More precisely,
only the $2$-hole-types such that the total number of
edges bounding the special holes is at most  $2k_1+2k_2+2$
have to be considered.

\section{Proof of equation (\ref{eq:i})}\label{sec:eq-i}

By equation \prettyref{eq:dZst-twotwo} we can rewrite \prettyref{eq:i}
as
\begin{equation*}
-2\sum_{\Gamma\in{\mathcal S}_z^{\leq 1}}
\Coeff_{z}^{-1}\langle\!\langle{\Gamma}\rangle\!
\rangle_{z\oplus\Lambda}^{[1]}
=
\gint{
{\xy\rgvertex[{\!}]{3}\loose1\loose2\loose3%
,(0,1.7)*+[F]{\Lambda^{-1}};(0,.5)**\dir{-}
\endxy}}\,,
\end{equation*}
where ${\mathcal S}_z^{\leq 1}$ denotes the set of hole types
such that the \(z\)-decorated hole is bounded by at most one edge.
There is only one such hole type. We compute
\begin{align*}
-2\Coeff_z^{-1}\left(\negquad
{\xy\rgvertex{11}\fence3{i}{i}\join7{10},(0,-.5)*{z}\endxy}
\negquad
\right)&=
-2\Coeff_z^{-1}\left(\frac{\U}{z+\Lambda_i}\right)\\
&=-2\U=
-2\U\negquad\negquad
{\xy,(0,-.3)*{}\rgvertex[\!]{3}\fence1{i}{i}\endxy}\negquad\,,
\end{align*}
so that
\begin{equation}\label{eq:det0}
-2\sum_{\Gamma\in{\mathcal S}_z^{\leq 1}}
\Coeff_{z}^{-1}\langle\!\langle{\Gamma}\rangle\!
\rangle_{z\oplus\Lambda}^{[1]}=
-2\U\gint{{\negquad\negquad\negquad
\xy,(0,-.3)*{}\rgvertex[\!]{3}\loose1\endxy
\negquad\negquad\negquad}}\,.
\end{equation}
Since the leg outgoing from the univalent vertex on the right-hand
side of equation \prettyref{eq:det0} must end in a trivalent vertex,
we find
\begin{equation*}
-2\!\!\sum_{\Gamma\in{\mathcal S}_z^{\leq 1}}
\Coeff_{z}^{-1}\langle\!\langle{\Gamma}\rangle\!
\rangle_{z\oplus\Lambda}^{[1]}=
-2\U\gint{{\xy\rgvertex{3}\loose1\loose2\loose3%
,(0,.7)*{\rgvertex[{\!}]{3}};(0,.5)**\dir{-}%
\endxy}}
\!\!=
\gint{{\xy\rgvertex[\!]{3}\loose1\loose2\loose3%
,(0,.7)*{\rgvertex[{\!}]{3}};(0,.5)**\dir{-}%
\endxy}}
\end{equation*}
On the other hand,
\begin{equation*}
{\xy\rgvertex[{\!}]{3}\loose1\fence2{i}{j}\fence3{j}{i}%
,(0,1.7)*+[F]{\Lambda^{-1}};(0,.5)**\dir{-}
\endxy}
={\Lambda_i}^{-1}={\xy\rgvertex[\!]{3}\loose1\fence2{i}{j}\fence3{j}{i}%
,(0,.7)*{\rgvertex[{\!}]{3}};(0,.5)**\dir{-}%
\endxy}
\end{equation*}
so
\begin{equation*}
\gint{{\xy\rgvertex[{\!}]{3}\loose1\loose2\loose3%
,(0,1.7)*+[F]{\Lambda^{-1}};(0,.5)**\dir{-}
\endxy}}
=\gint{{\xy\rgvertex[\!]{3}\loose1\loose2\loose3%
,(0,.7)*{\rgvertex[{\!}]{3}};(0,.5)**\dir{-}%
\endxy}}
\end{equation*}
and equation \prettyref{eq:i} follows.

\section{Proof of equation (\ref{eq:ii})}\label{sec:eq-ii}

Th proof of equation \prettyref{eq:ii} goes along the same
lines of the proof of equation, \prettyref{eq:i}, but it is
a bit more involved. We begin by rewriting
\prettyref{eq:ii} by means of the Witten's formulas
(\ref{eq:dZst-twotwo}-\ref{eq:dZst-twotwoo}) for the derivatives of the
partition function:
\begin{align}
\notag
&-14\sum_{\Gamma\in{\mathcal S}_z^{\leq 7}}
  \Coeff_{z}^{-7}\langle\!\langle{\Gamma}\rangle\!
\rangle_{z\oplus\Lambda}^{[1]}
-6
\sum_{\Gamma\in{\mathcal S}_{z,w}^{\leq4}}\Coeff_{w}^{-1}
\Coeff_{z}^{-3}
  \langle\!\langle{\Gamma}\rangle\!
\rangle_{w\oplus z\oplus\Lambda}^{[1,1]}-\\
\notag&\qquad
-6\tr\Lambda
\sum_{\Gamma\in{\mathcal S}_z^{\leq 3}}
  \Coeff_{z}^{-3}\langle\!\langle{\Gamma}
\rangle\!\rangle_{z\oplus\Lambda}^{[1]}
-2\tr\Lambda^3
\sum_{\Gamma\in{\mathcal S}_z^{\leq1}}
  \Coeff_{z}^{-1}\langle\!\langle{\Gamma}
\rangle\!\rangle_{z\oplus\Lambda}^{[1]}
-\\
\label{eq:final-ii}
&\qquad-
\gint{
{\xy\rgvertex[{\!}]{3}\loose1\loose2\loose3%
,(0,1.7)*+[F]{\Lambda^5};(0,.5)**\dir{-}
\endxy}
}+\bigl(2\tr\Lambda^{4}-2\tr\Lambda^3\tr\Lambda+(\tr\Lambda^2)^2
\bigr)\cdot\langle\!\langle\emptyset\rangle\!\rangle_\Lambda^{}=0\,,
\end{align}
where ${\mathcal S}_z^{\leq l}$ denotes the set of isomorphism classes
of hole
types such that the $z$-decorated hole is bounded by at most $l$ edges, and
${\mathcal S}_{z,w}^{\leq 4}$ denotes the set of isomorphism classes of
$2$-holes
types such that the total number of edges bounding the $z$- and
$w$-decorated
holes is at most four.

 By the Feynman rules in the 't~Hooft-Kontsevich model, the Laurent
coefficients of the amplitude of an hole type are
polynomials in the variables $\Lambda_i$'s. Indeed:
\begin{enumerate}
\item each ordinary vertex brings a
factor \(\U/2)\);
\item each internal edge bordering the \(z\)-decorated hole on
both sides contributes a factor \(1/z\);
\item the other internal edges  contribute factors
of the form \(2/(z+\Lambda_i)\) for \(i=1,\dots,N\).
\end{enumerate}
Similarly, also the Laurent coefficients of the amplitudes
of a two-holes type are polynomials in the variables
$\Lambda_i$'s.

A convenient way to represent graphically  these Laurent
coefficients is to enlarge the class of special vertices by
adding special vertices decorated by polynomials. This is
formally done as follows.

Let \(\varphi\) be a polynomial in
\(\setC[\theta_1,\theta_2,\dots,\theta_n]\); we say that the
polynomial \(\varphi\) is cyclically invariant if it is
invariant with respect to the natural action of the cyclic group
\({\mathbb
  Z}/n{\mathbb Z}\) on the coordinates. We will consider an
\(n\)-valent special vertex decorated by the polynomial
\(\varphi\):
\begin{equation*}
  \xyc
  \rgvertex[\varphi]{5}\genericfences
  \endxyc
  \mapsto
  \varphi(\Lambda_{i_1},\Lambda_{i_2},\dots,\Lambda_{i_n})
\end{equation*}
This Feynman rule is well defined due to the cyclical invariance of
\(\varphi\).

We also consider special vertices decorated by non-cyclically invariant
polynomials:
\begin{equation*}
  \xyc\rgvertex[\varphi]{5}\genericfences\cilia{2}{3}\endxyc
  \mapsto
  \varphi(\Lambda_{i_1},\Lambda_{i_2},\dots,\Lambda_{i_n}),
\end{equation*}
 The r{\^o}le of the ``$\bigstar$'' mark is precisely to break the
cyclical symmetry of the graphical element,
so that the ``$\bigstar$'' tells which indeterminate ---among those
corresponding to indices decorating holes around the vertex--- comes
first. Note that, when dealing with $\bigstar$-marked special vertices,
there
is actually
no need that the valence $n$ of the vertex equals the the number $\nu$ of
variables of the polynomial decorating it. Indeed, if $\nu < n$ then the
above
Feynman rule still makes sense; if $\nu > n$ then wrap around the
vertex as many times as needed.  

Using these notations, the Laurent coefficients we are interested in
are easily written in the form
\begin{equation*}
\sum_{\Xi}p_{\Xi}(\tr\Lambda^*)Z_\Lambda(\Xi)
\end{equation*}
where the $\Xi$'s are disjoint unions  
of special vertices decorated by polynomials,
and the $p_\Xi^{}$'s are polynomials in the \emph{positive} traces of
$\Lambda$;
see
\cite{fiorenza-murri;matrix-integrals} for several computation
examples; we call $\Xi$
a
\emph{cluster} of special vertices.

Also the
term
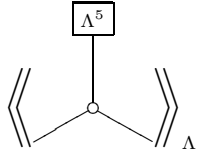
\begin{equation*}
\gint{
{\xy\rgvertex[{\!}]{3}\loose1\loose2\loose3%
,(0,1.7)*+[F]{\Lambda^5};(0,.5)**\dir{-}
\endxy}
}
\end{equation*}
appearing in \prettyref{eq:final-ii} can be written as the
expectation value of a special vertex decorated by a polynomial.
Indeed,
\begin{equation*}
{\xy\rgvertex[{\!}]{3}\fence3{i}{j}\fence2{j}{i}\loose1%
,(0,1.7)*+[F]{\Lambda^5};(0,.5)**\dir{-}
\endxy}={\Lambda_i}^5=
{\xy\rgvertex[{\theta_1}^5]{6}\fence1{i}{j}\fence4{j}{i}
\cilia{1}{2}
\endxy}
\end{equation*}
so that
\begin{equation*}
\gint{
{\xy\rgvertex[{\!}]{3}\loose1\loose2\loose3%
,(0,1.7)*+[F]{\Lambda^5};(0,.5)**\dir{-}
\endxy}
}=
\gint{
{\xy\rgvertex[{\theta_1}^5]{6}\loose1\loose4
\cilia{1}{2}
\endxy}
}
\end{equation*}
Therefore, equation \prettyref{eq:final-ii} is reduced to showing
that
\begin{equation}\label{eq:canonical-form-eq}
\sum_{\Xi}q^{}_\Xi(\tr\Lambda^*)\langle\!\langle{\Xi}
\rangle\!\rangle^{}_\Lambda=0
\end{equation}
where the $\Xi$ are suitable clusters of special vertices and the
$q^{}_\Xi$ are polynomials in the positive traces of $\Lambda$. We
remark that all the clusters
$\Xi$ and the polynomials $q^{}_\Xi$ are explicitly computable: one
just has to write down all the hole types and the two-holes types
involved in equation \prettyref{eq:final-ii}.

To prove equation \prettyref{eq:canonical-form-eq} we will
use the Feynman diagrammatic manipulations introduced in
\cite{fiorenza-murri;matrix-integrals}. These manipulations will
reduce the proof to a completely straightforward, but
long and tedious, computation. Therefore, after describing the
general theory, we will explicitly
compute only a few terms of \prettyref{eq:canonical-form-eq},
leaving the remaining computations as an exercise to the reader.

 We define
the \emph{total degree} of a term of the form
\begin{equation*}
(\tr\Lambda^0)^{d_0}
(\tr\Lambda)^{d_1}\cdots (\tr\Lambda^n)^{d_n}\langle\!\langle{\Xi}
\rangle\!\rangle^{}_\Lambda
\end{equation*}
as the integer $\deg\Xi+\sum_{i} i\,d_i$, where $\deg\Xi$ is the
sum of the degrees of the polynomials decorating the vertices of
$\Xi$.

The key idea (due to Witten \cite{witten;kontsevich-model}) is that one
can lower the total degrees of the terms appearing in equation
\prettyref{eq:canonical-form-eq} by the following remark. For any
$\psi
\in
\setC[\theta_1,
\ldots,
\theta_n]$, let
  \(u_\psi \in \setC[\theta_1,\dots,\theta_n]\) be the
polynomial
\begin{equation*}
   u_\psi(\theta_1,\theta_2,\dots,\theta_n):=
    (\theta_n+\theta_1)\cdot\psi(\theta_1,\theta_2,\dots,\theta_n),
  \end{equation*}  
Then
\begin{equation*}
 {\xyc\rgvertex[u_\psi]{7}%
 \fence1{i_n}{i_1}\fence2{i_1}{i_2}%
\fence3{i_2}{i_3}\fence4{i_4}{\dots}\missing5\fence7{\dots}{i_n}%
        \cilia{1}{2}\endxyc}\mapsto
(\Lambda_{i_1}+\Lambda_{i_n})\psi(\Lambda_1,\dots,\Lambda_{i_n})
  \end{equation*}
and
 then we can use the factor $(\Lambda_{i_1}+\Lambda_{i_n})$ to
cancel the factor $2/(\Lambda_{i_1}+\Lambda_{i_n})$ coming from the
edge stemming from the vertex
  just before the ciliation (in the cyclic order of the vertex).
Now, assume we are dealing with the expectation value
\begin{equation*}
\gint{{\xyc\rgvertex[u_\psi]{7}%
        \loose1\loose2\loose3\loose4\missing5\loose7%
        \cilia{1}{2}\endxyc}}
\end{equation*}
  By definition, this expectation value  is a sum
  over ribbon graphs (with a distinguished sub-diagram); for any
\(\Gamma\)
  in this sum, the edge stemming from the vertex
  just before the ciliation must
  \emph{either} end at another ---distinct--- vertex \emph{or} make a
  loop. Therefore, using the definition of expectation value again,
  \begin{multline}\label{eq:expansion}
    \gint{{\xyc\rgvertex[u_\psi]{13}\loose1%
        \loose3%
        \loose5%
        \loose{10}%
        \loose{12}%
        \missing7%
        \cilia{1}{3}%
        \endxyc}}=\gint{{
        \xyc
        \rgvertex[u_\psi]{13}%
        \loose1%
        \loose3%
        \loose5%
        \loose{10}%
        \loose{12}%
        \missing7%
        \cilia{1}{3}%
        ,(1.3,0);(1.6,.5)**\dir{-},(1.3,0);(1.6,-.5)**\dir{-}%
        ,(1.3,0)*{\bullet}%
        \endxyc}}
     +\gint{{\xyc
        \rgvertex[u_\psi]{13}%
        \join{1}{3}
        \loose5%
        \loose{10}%
        \loose{12}%
        \missing7%
        \cilia{1}{3}%
        \endxyc}}+
    \\
    +\sum_{j=1}^{n-3}\gint{{
        \xyc
        \rgvertex[u_\psi]{13}%
        \loose2%
        \loose4%
        \loose6%
        \loose{10}%
        \missing{3}%
        \loose{12}%
        \missing7%
        \cilia{1}{2}%
        ,(1.5,1.5)*{j\ {\rm edges,}}
        ,(1.5,1.2)*{0<j<n-2}
        ,(.405,0);(.85,.35)**\crv{(1.2,-0.2)}%
        ,(.75,.48);(.15,.92)**\crv{(.6,.7)}
        ,(-.15,.39);(.05,.95)**\crv{(-.3,1.1)}
        \endxyc
      }}
   +\gint{{\xyc
        \rgvertex[u_\psi]{13}%
        \join{1}{12}
        \loose5%
        \loose{10}%
        \loose{3}%
        \missing7%
        \cilia{1}{3}%
        \endxyc}}
  \end{multline}
  Now, each of the terms at right-hand side of \eqref{eq:expansion}
  above, can be rewritten as the expectation
  value of a linear combination (over $\setC[\tr
\Lambda^*]$) of
  clusters; indeed, one can directly compute:
  \begin{align}
    \label{eq:C1}\tag{C1}
    \xyc
    \rgvertex[u_\psi]{13}%
    \loose1%
    \loose3%
    \loose5%
    \loose{10}%
    \loose{12}%
    \missing7%
    \cilia{1}{3}%
    ,(1.3,0);(1.6,.5)**\dir{-},(1.3,0);(1.6,-.5)**\dir{-}%
    ,(1.3,0)*{\bullet}
    \endxyc
    &=
    \U \cdot
    \underbrace{
      \xyc
      \rgvertex[\psi]{13}%
      \loose2%
      \loose{13}%
      \loose3%
      \loose5%
      \loose{10}%
      \loose{12}%
      \missing7
    \cilia{2}{3}
      \endxyc
    }_{\text{$(n+1)$-valent}}\qquad;
    \\ \label{eq:C2}\tag{C2}
    \xyc
    \rgvertex[u_\psi]{13}%
    \join{1}{3}
    \loose5%
    \loose{10}%
    \loose{12}%
    \missing7%
    \cilia{1}{3}%
    \endxyc
    &=
    2\sum_{h=0}^k\tr\Lambda^h\cdot
    \underbrace{
      \xyc
      \rgvertex[\psi_h]{13}%
      \loose5%
      \loose{10}%
      \loose{12}%
      \missing7%
      ,(.5,.5)*{\scriptstyle{\bigstar}}
      \endxyc
    }_{\text{$(n-2)$-valent}}
    ,
\end{align}
where the polynomials $\psi_h$ are defined by the equation
\begin{equation*}
\psi(\theta_1, \dots, \theta_{n-1}, \theta_2)=
\sum_{h=0}^k {\theta_1}^h \psi_h(\theta_2, \dots,
\theta_{n-1});
\end{equation*}
\begin{equation}
    \label{eq:C3}\tag{C3}
    {\xyc
      \rgvertex[u_\psi]{13}%
      \loose2%
      \loose4%
      \loose6%
      \loose{10}%
      \missing{3}%
      \loose{12}%
      \missing7%
      \cilia{1}{2}%
      ,(1.5,1.5)*{j\ {\rm edges,}}
      ,(1.5,1.2)*{0<j<n-2}
      ,(.405,0);(.85,.35)**\crv{(1.2,-0.2)}%
      ,(.75,.48);(.15,.92)**\crv{(.6,.7)}
      ,(-.15,.39);(.05,.95)**\crv{(-.3,1.1)}
      \endxyc}
    =
    2\sum_{h=0}^k
    \underbrace{
      \xyc
      \rgvertex[\phi_{h}']{13}%
      \loose2%
      \loose4%
      \missing3%
      \ciliamedio{1}{2}%
      \endxyc
    }_{\text{$j$-valent}}
    \ \underbrace{
      \xyc
      \rgvertex[\phi_{h}'']{13}%
      \loose6%
      \loose{12}\loose{10}%
      \missing7%
      \ciliamedio{1}{2}
      \endxyc
    }_{\text{$(n-j-2)$-valent}}
    ,
  \end{equation}
where the polynomials $\phi_{h}'$, $\phi_{h}''$ are defined by
\begin{multline*}
\psi(\theta_1,
    \dots, \theta_j, \theta_1, \theta_{j+2}, \dots, \theta_{n-1},
    \theta_{j+2}) =\\
=\sum_{h=0}^k \phi_{h}'(\theta_1, \dots, \theta_j) \cdot
    \phi_{h}''(\theta_{j+2}, \dots, \theta_{n-1})
 ;
\end{multline*}
  \begin{equation}
    \label{eq:C4}\tag{C4}
    \xyc
    \rgvertex[u_\psi]{13}%
    \join{1}{12}
    \loose5%
    \loose{10}%
    \loose{3}%
    \missing7%
    \cilia{1}{3}%
    \endxyc
    =
    2\sum_{h=0}^k\tr\Lambda^h\cdot
    \underbrace{
      \xyc
      \rgvertex[\eta_h]{13}%
      \loose5%
      \loose{10}%
      \loose{3}%
      \missing7%
      ,(.5,-.5)*{\scriptstyle{\bigstar}}
      \endxyc
    }_{\text{$(n-2)$-valent}}
    ,
  \end{equation}
 where the polynomials
  $\eta_h$ are
  defined by:
  \begin{equation*}
    \psi(\theta_1, \dots, \theta_{n-2}, \theta_1, \theta_n)=
    \sum_{h=0}^k {\theta_n}^h \eta_h(\theta_1, \dots, \theta_{n-2}).
  \end{equation*}
 
  The same argument also works in the more general case of a
cluster of special vertices such that one of the polynomials
decorating the vertices has the form $u_\psi$. In this case,
  a new graphical element may appear,
which
  is not listed in equations \eqref{eq:C1}--\eqref{eq:C4} above;
  namely, that the ciliated edge connects the chosen vertex to another
  one in the same cluster. Direct computation again gives:
  \begin{equation*}
    \label{eq:C1'}\tag{C1'}
    \xyc
    \rgvertex[u_\psi]{13}%
    \loose1%
    \loose3%
    \loose5%
    \loose{10}%
    \loose{12}%
    \missing7%
    \cilia{1}{3}%
    ,(1.6,0)*{\ \zeta\ }*\cir{};(2.15,.8)**\dir{-}%
    ,(1.6,0)*{\ \zeta\ }*\cir{};(2.15,-.8)**\dir{-}%
    ,(2.1,0)*{\scriptstyle{\cdots}}
    ,(3.2,0)*{\Bigg\}\text{$j$ edges}}
    \endxyc
    = 2\cdot\underbrace{
      \xyc
      \rgvertex[\psi*\zeta]{13}%
      \loose2%
      \loose{13}%
      \loose3%
      \loose5%
      \loose{10}%
      \loose{12}%
      \ciliamedio{2}{3},(-1,.2)*{\cdots},(1,0)*{\cdots}
      \endxyc
    }_{\text{$(n+j-1)$-valent}}
    ,
  \end{equation*}
  where the polynomial $\psi*\zeta$ is defined by
  \begin{equation*}
    (\psi*\zeta)(\theta_1,\dots,\theta_{n+j})=
    \psi(\theta_1,\dots,\theta_n)\cdot\zeta(\theta_n,\theta_{n+1},
    \dots,\theta_{n+j},\theta_1).
  \end{equation*}
 
By this edge-contracting procedure we can lower the total degree of
expressions of the form
\begin{equation*}
p_\Xi\langle\!\langle\Xi\rangle\!\rangle^{}_\Lambda
\end{equation*}
as soon as one of the vertices of $\Xi$ is decorated by a
polynomial of the form $u_\psi$. In order to apply it to equation
\prettyref{eq:canonical-form-eq} we have to change the polynomials
decorating the vertices in \prettyref{eq:canonical-form-eq} with
polynomials of the form
$u_\psi$. This can be done by using the following trick: if
$\varphi_1^{}(\theta_1,\dots,\theta_n)$ and
$\varphi_2^{}(\theta_1,\dots,\theta_n)$ are two polynomials such
that
\begin{equation*}
\sum_{\sigma\in\setZ/n\setZ}
\varphi_1^{}(\theta_{\sigma(1)},\dots,\theta_{\sigma(n)})
=
\sum_{\sigma\in\setZ/n\setZ}
\varphi_2^{}(\theta_{\sigma(1)},\dots,\theta_{\sigma(n)})
\end{equation*}
then
\begin{equation*}
\gint{{\xyc\rgvertex[\varphi_1^{}]{7}%
        \loose1\loose2\loose3\loose4\missing5\loose7%
        \cilia{1}{2}\endxyc}}
=
\gint{{\xyc\rgvertex[\varphi_2^{}]{7}%
        \loose1\loose2\loose3\loose4\missing5\loose7%
        \cilia{1}{2}\endxyc}}
\end{equation*}
Indeed, if we set
\begin{equation*}
\overline{\varphi}(\theta_1,\dots,\theta_n):=
\sum_{\sigma\in\setZ/n\setZ}
\varphi_1^{}(\theta_{\sigma(1)},\dots,\theta_{\sigma(n)})
\end{equation*}
then
\begin{equation*}
\gint{{\xyc\rgvertex[\varphi_1^{}]{7}%
        \loose1\loose2\loose3\loose4\missing5\loose7%
        \cilia{1}{2}\endxyc}}
=
\gint{{\xyc\rgvertex[\overline{\varphi}]{7}%
        \loose1\loose2\loose3\loose4\missing5\loose7%
        \endxyc}}
\end{equation*}
and the same equality holds with $\varphi_2^{}$ in place of
$\varphi^{}_1$.
Therefore, a way of changing a polynomial $\varphi$ decorating a
vertex into a polynomial of the form $u_\psi$ is the following:
\begin{enumerate}
\item make a cyclically invariant polynomial out of $\varphi$ by
setting
$$\overline{\varphi}^{}(\theta_1,\dots,\theta_n):=
\sum_{\sigma\in\setZ/n\setZ}
\varphi^{}(\theta_{\sigma(1)},\dots,\theta_{\sigma(n)});$$
\item find a cyclic decomposition
\begin{equation}\label{eq:cyc-dec}
    \overline{\varphi}(\theta_1,\dots,\theta_n) =
 \sum_{\sigma\in\setZ/n\setZ} u_\psi(
    \theta_{\sigma(1)}, \dots, \theta_{\sigma(n)}).
  \end{equation}
\end{enumerate}
Clearly, not every cyclically invariant polynomial admits a cyclic
decomposition of the form \prettyref{eq:cyc-dec}; for instance the
polynomial
$\overline{\varphi}(\theta_1,\theta_2)={\theta_1}^2+{\theta_2}^2$
has no such decomposition. Therefore, in general one could be not
able to change the polynomial $\varphi$ into a polynomial of the
form
$u_\psi$. However, explicit computations show that all the
polynomials occurring in the proof of equation
\prettyref{eq:final-ii} can actually be changed in the form $u_\psi$ by
using the trick described above. We address the reader interested
in dealing with cyclic polynomials which do not admit a cyclic
decomposition of the form \prettyref{eq:cyc-dec} to
\cite{fiorenza-murri;matrix-integrals}.

We now work out explicitly the first steps in the proof of
equation \prettyref{eq:final-ii} to show how the above general theory
applies.

 All term in equation
\prettyref{eq:canonical-form-eq} have at most total degree equal to
$6$. Moreover, there is only one term of total degree $6$: it is
the term coming form the unique hole type such that the
$z$-decorated hole is bounded by exactly one edge.
We have
\begin{align*}
-14\Coeff_z^{-7}\left(
{\xy\rgvertex{11}\fence3{i}{i}\join7{10},(0,-.5)*{z}\endxy}\right)&=
-14\Coeff_z^{-7}\frac{\sqrt{-1}}{z+\Lambda_i}=
\\
&=-14\sqrt{-1}\,\,{\Lambda_i}^6=
-14{\U}\cdot\negquad\negquad
{\xy\rgvertex[{\theta_1}^{6}]{3}\fence1{i}{i}\endxy}\negquad
\end{align*}
So, the degree $6$ term in
\prettyref{eq:canonical-form-eq} is
\begin{equation*}
-14\U\gint{\negquad\negquad
{\xy\rgvertex[{\theta_1}^{6}]{3}\loose1\endxy}\negquad\negquad}
\end{equation*}
 The
polynomial
${\theta_1}^6$ is trivially cyclically
decomposable:
$$\displaystyle{{\theta_1}^6=(\theta_1+\theta_1)\cdot
\frac{{\theta_1}^5}{2}}.$$ Therefore, by the
edge-contraction procedure described above, we
find
\begin{equation*}
-14\U\gint{\negquad\negquad
{\xy\rgvertex[{\theta_1}^{6}]{3}\loose1\endxy}\negquad\negquad}=
7\gint{
{\xy\rgvertex[{\theta_1}^5]{6}\loose1\loose4
\cilia{1}{2}\endxy}
}
\end{equation*}
This way we have reduced the degree $6$ term to a degree $5$
term. In equation \prettyref{eq:canonical-form-eq} there are
other two degree $5$ terms. One is the term coming from the
unique hole type such that the $z$-decorated hole is bounded
by exactly two edges. We have
\begin{align*}
-14\Coeff_{z}^{-7}\left({\xy
,(0,-.8)*{\bullet};(0,.8)*{\bullet}**\crv{(-1.3,0)}
,(0,-.8)*{\bullet};(0,.8)*{\bullet}**\crv{(1.3,0)}
,(0,-.8)*{\bullet};(0,-1.7)**\dir{-}
,(0,.8)*{\bullet};(0,1.7)**\dir{-}
,(0,0)*{z},(-.2,-1.6)*{\scriptstyle{i_1}}
,(-.2,1.6)*{\scriptstyle{i_1}}
,(.2,-1.6)*{\scriptstyle{i_2}}
,(.2,1.6)*{\scriptstyle{i_2}}
\endxy}\right)&=
14\Coeff_z^{-7}\frac{1}{(z+\Lambda_{i_1})(z+\Lambda_{i_2})}\\
&\negquad
=-14\displaystyle{\sum_{n_1+n_2=5}}{\Lambda_{i_1}}^{n_1}
{\Lambda_{i_2}}^{n_2}
=
{\xy\rgvertex[\eta^{}_{2}]{6}\fence1{i_1}{i_2}
\fence4{i_2}{i_1}
\endxy}
\end{align*}
where
\begin{equation*}
\eta^{}_2(\theta_1,\theta_2)=-14\bigl({\theta_1}^5+
{\theta_1}^4{\theta_2}+{\theta_1}^3{\theta_2}^2+
{\theta_1}^2{\theta_2}^3+
{\theta_1}{\theta_2}^4+{\theta_2}^5\bigr)\,.
\end{equation*}
The other degree $5$ term is the one coming from
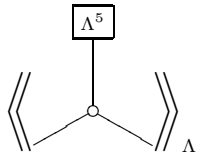
\begin{equation*}
\gint{
{\xy\rgvertex[{\!}]{3}\loose1\loose2\loose3%
,(0,1.7)*+[F]{\Lambda^5};(0,.5)**\dir{-}
\endxy}
}
\end{equation*}
Indeed, we have seen at the beginning of this section that this
expectation value equals
\begin{equation*}
\gint{
{\xy\rgvertex[{\theta_1}^5]{6}\loose1\loose4
\cilia{1}{2}\endxy}
}
\end{equation*}

Summing up, we obtain that the contribution from terms
of total degree at least
$5$ in equation \prettyref{eq:canonical-form-eq} is reduced
to:
\begin{equation*}
\gint{
{\xy\rgvertex[\overline{\varphi}_2]{6}\loose1\loose4
\endxy}
}
\end{equation*}
where
\begin{equation*}
\overline{\varphi}_2^{}(\theta_1,\theta_2)=-8{\theta_1}^5-14
{\theta_1}^4{\theta_2}-14{\theta_1}^3{\theta_2}^2
-14{\theta_1}^2{\theta_2}^3-14{\theta_1}{\theta_2}^4
-8{\theta_2}^5.
\end{equation*}
 The polynomial $\varphi_2^{}$ is
cyclically decomposable: we have\footnote{
The computer program Maple V has been uses for this and the
following computations.}
\begin{equation*}
\overline{\varphi}_2^{}(\theta_1,\theta_2)=
\sum_{\sigma\in{\setZ}/2\setZ}u_{\psi_2}(
\theta_{\sigma(1)},\theta_{\sigma(2)})\,,
\end{equation*}
where
\begin{equation*}
\psi_2^{}(\theta_1,\theta_2)=-4{\theta_1}^4-
3{\theta_1}^3{\theta_2}-4{\theta_1}^2{\theta_2}^2-
3{\theta_1}{\theta_2}^3-4{\theta_2}^4\,.
\end{equation*}

Applying the edge-contraction procedure again, we get
\begin{align*}
\gint{
{\xy\rgvertex[\overline{\varphi}_2]{6}\loose1\loose4
\endxy}
}&=
\U
\gint{
{\xy\rgvertex[\psi^{}_2]{3}\loose1\loose2\loose3\cilialontano{1}{2}
\endxy}
}+\\
&\phantom{mmmmmmm}+
\left(-16\tr\Lambda^4-12\tr\Lambda^3\tr\Lambda-
8(\tr\Lambda^2)^2\right)
\langle\!\langle\emptyset\rangle\!\rangle_\Lambda^{}\\
&\negquad\negquad\negquad\negquad\negquad\negquad=
\U
\gint{
{\xy\rgvertex[\overline{\psi}^{}_2]{3}\loose1\loose2\loose3\endxy}
}
-
\left(16\tr\Lambda^4+12\tr\Lambda^3\tr\Lambda+
8(\tr\Lambda^2)^2\right)
\langle\!\langle\emptyset\rangle\!\rangle_\Lambda^{}
\end{align*}
where
\begin{align*}
\overline{\psi}^{}_2(\theta_1,\theta_2,\theta_3)&=
\psi^{}_2(\theta_1,\theta_2)+\psi^{}_2(\theta_2,\theta_3)
+\psi^{}_2(\theta_2,\theta_1)\\
&=
-8{\theta_1}^4-
3{\theta_1}^3{\theta_2}-4{\theta_1}^2{\theta_2}^2-
3{\theta_1}{\theta_2}^3-8{\theta_2}^4
-3{\theta_2}^3{\theta_3}-\\
&\qquad -4{\theta_2}^2{\theta_3}^2-
3{\theta_2}{\theta_3}^3
-8{\theta_3}^4-
3{\theta_3}^3{\theta_1}-4{\theta_3}^2{\theta_1}^2-
3{\theta_3}{\theta_1}^3\,.
\end{align*}
This way we have reduced the total contribution coming from
terms of total degree at least $5$ to a sum of
terms of total degree $4$. There are other terms of
total degree $4$ in equation
\prettyref{eq:canonical-form-eq}. One is the term
\begin{equation*}
\left(2\tr\Lambda^{4}-2\tr\Lambda^3\tr\Lambda+(\tr\Lambda^2)^2
\right)\cdot
\langle\!\langle\emptyset\rangle\!\rangle_\Lambda^{}
\end{equation*}
appearing in equation \prettyref{eq:final-ii}. The others are the
terms coming from the hole types whose $z$-decorated hole is
bounded by exactly $3$ edges. There are three such hole types. The
first of them gives
\begin{align*}
-14\Coeff_z^{-7}\left(\!
{\xy\rghole[z]{3}\fence1{i_1}{i_2}\fence2{i_2}{i_3}
\fence3{i_3}{i_1}\endxy}\!\right)&=
14\Coeff_z^{-7}\frac{\U}
{(z+\Lambda_{i_1})(z+\Lambda_{i_2})(z+\Lambda_{i_3})}\\
&=
14\U
\displaystyle{\sum_{n_1+n_2+n_3=4}}{\Lambda_{i_1}}^{n_1}
{\Lambda_{i_2}}^{n_2}{\Lambda_{i_3}}^{n_3}\\
&=
{\xy\rgvertex[\eta^{}_{3}]{3}\fence1{i_1}{i_2}
\fence2{i_2}{i_3}\fence3{i_3}{i_1}
\endxy}
\end{align*}
where
$\displaystyle{\eta^{}_{3}(\theta_1,\theta_2,\theta_3)=14\U
\sum_{n_1+n_2+n_3=4}{\theta_{1}}^{n_1}{\theta_{2}}^{n_2}
{\theta_{3}}^{n_3}}$.

The second hole type with three edges bounding the $z$-decorated
hole gives no contribution to total degree $4$ terms. Indeed,
\begin{equation*}
-14\Coeff_{z}^{-7}\left({\xyc/r.5cm/:
      \vloop~{(-1,2)}{(1,2)}{(-1,1)}{(1,1)}%
      |{*{}="m"},"m"*{\bullet};(0,1)**\dir{-},%
      (-1,1);(-1,-1)**\dir{-},%
      \vcross~{(0,1)}{(1,1)}{(0,-1)}{(1,-1)},%
      \vloop~{(-1,-2)}{(1,-2)}{(-1,-1)}{(1,-1)}%
      |{*{}="m"},"m"*{\bullet};(0,-1)**
      \dir{-},(-0.2,0)*{z}%
      \endxyc}\right)=\frac{7}{2}\Coeff_z^{-7}\frac{1}{z^3}=0
\phantom{mmmmmmmmmmm}
\end{equation*}

Finally, the third hole type gives
\begin{align*}
&-14\Coeff_z^{-7}\left({\xy,(0,-.2)*!C\xybox{
        ,(3.1,1.4);(3.05,1.2)**\crv{(3.13,1.5)&(3.5,1.866)%
          &(4,2)&(4.5,1.866)&(5,1)&(4,0)&(3,0.7)}
        ,(3.068,1.0);(2.5,1.4)**\crv{(3.2,1.07)&(3.3,1.09)%
          &(3.3,1.15)&(3.07,1.35)}
        ,(1.9,1.4);(1.95,1.2)**\crv{(1.87,1.5)&(1.5,1.866)%
          &(1,2)&(.5,1.866)&(0,1)&(1,0)&(2,0.7)}
        ,(1.932,1.0);(2.5,1.4)**\crv{(1.8,1.07)&(1.7,1.09)%
          &(1.7,1.15)&(1.93,1.35)}%
     ,(3.064,1.0)*{\bullet},(1.938,1.0)*{\bullet}%
,(1.15,1.2)*{z}%
       }\endxy}\right)\\
&\qquad-14\sum_{i_1,i_2}\Coeff_z^{-7}\left({\xy,(0,-.2)*!C\xybox{
        ,(3.1,1.4);(3.05,1.2)**\crv{(3.13,1.5)&(3.5,1.866)%
          &(4,2)&(4.5,1.866)&(5,1)&(4,0)&(3,0.7)}
        ,(3.068,1.0);(2.5,1.4)**\crv{(3.2,1.07)&(3.3,1.09)%
          &(3.3,1.15)&(3.07,1.35)}
        ,(1.9,1.4);(1.95,1.2)**\crv{(1.87,1.5)&(1.5,1.866)%
          &(1,2)&(.5,1.866)&(0,1)&(1,0)&(2,0.7)}
        ,(1.932,1.0);(2.5,1.4)**\crv{(1.8,1.07)&(1.7,1.09)%
          &(1.7,1.15)&(1.93,1.35)}%
     ,(3.064,1.0)*{\bullet},(1.938,1.0)*{\bullet}%
,(1.15,1.2)*{z},(1.1,.1)*{\scriptstyle{i_1}}%
,(4.1,.1)*{\scriptstyle{i_2}}%
       }\endxy}\right)\\
\\
&\qquad
=14\Coeff_z^{-7}\sum_{i_1,i_2}\frac{1}{z(z+\Lambda_{i_1})
(z+\Lambda_{i_1})}\\
&\\
&\qquad=
2\biggl(14\tr\Lambda^0\tr\Lambda^4+14\tr\Lambda\tr\Lambda^3+
7(\tr\Lambda^2)^2\biggr)
\langle\!\langle\emptyset\rangle\!\rangle_\Lambda^{}
\end{align*}
Summing up (taking into account the cardinalities of the automorphism groups
involved) the total contribution from terms of total degree at least $4$
is just
\begin{equation*}
\gint{\xy\rgvertex[\overline{\varphi}^{}_{4}]{3}\loose1\loose2\loose3\endxy}
\end{equation*}
where
\begin{align*}
\overline{\varphi}^{}_{4}(\theta_1,\theta_2,\theta_3)=\U\bigl(
&
6{\theta_1}^4+11{\theta_1}^3{\theta_2}
+10{\theta_1}^2{\theta_2}^2
+11{\theta_1}{\theta_2}^3+6{\theta_2}^4
+\\
&\!\!+11{\theta_2}^3{\theta_3}+10{\theta_2}^2{\theta_3}^2+
11{\theta_2}{\theta_3}^3
+6{\theta_3}^4+11{\theta_3}^3{\theta_1}+\\
&\!\!+10{\theta_3}^2{\theta_1}^2+11{\theta_3}{\theta_1}^3
+14\theta_1\theta_2{\theta_3}^2
+14\theta_2\theta_3{\theta_1}^2+\\
&\!\!
+14\theta_3\theta_1{\theta_2}^2
\bigr)
\end{align*}
One computes
$\overline{\varphi}_3^{}(\theta_1,\theta_2,\theta_3)=
\sum_{\sigma\in{\setZ}/3\setZ}u_{\psi_3}(
\theta_{\sigma(1)},\theta_{\sigma(2)},\theta_{\sigma(3)})$,
where
\begin{align*}
\psi_3^{}(\theta_1,\theta_2,\theta_3)&=
\frac{\U}{4}\bigl(16{\theta_2}^2\theta_3+30{\theta_2}{\theta_3}^2
+15{\theta_3}^3+8{\theta_1}{\theta_2}^2
-4{\theta_1}{\theta_2}{\theta_3}+\\
&\qquad+15{\theta_1}{\theta_3}^2
+10{\theta_1}^2{\theta_2}+4{\theta_2}^3+
{\theta_1}^2{\theta_32}+9{\theta_1}^3\bigr)\,.
\end{align*}
Moreover,
\begin{equation*}
\psi^{}_3(\theta_1,\theta_2,\theta_2)=
\frac{\U}{4}(65{\theta_2}^3
+19{\theta_1}{\theta_2}^2
+11{\theta_1}^2{\theta_2}+9{\theta_1}^3
)
\end{equation*}
\begin{equation*}
\psi^{}_3(\theta_1,\theta_1,\theta_3)=
\frac{\U}{4}(45{\theta_1}{\theta_3}^2
+13{\theta_1}^2
{\theta_3}
+
31{\theta_1}^3+15{\theta_3}^3)
\end{equation*}
So, applying the edge-contraction procedure again, we find
\begin{align*}
\gint{\xy\rgvertex[\overline{\varphi}^{}_{4}]{3}\loose1\loose2\loose3\endxy}
&=\U\gint{
\xy\rgvertex[\overline{\psi}^{}_{3}]{4}\loose1\loose2\loose3\loose4\endxy
}+12\U\tr\Lambda^3
\gint{\negquad\negquad\negquad
{\xy,(0,-.3)*{}%
\rgvertex[\!]{3}\loose1\endxy}\negquad\negquad\negquad}+\\
&\negquad+28\U\tr\Lambda^2
\gint{\negquad\negquad
{\xy\rgvertex[{\theta_1}]{3}\loose1\endxy}\negquad\negquad}
\!+16\U\tr\Lambda
\gint{\negquad\negquad
{\xy\rgvertex[{\theta_1}^{2}]{3}\loose1\endxy}\negquad\negquad}\!+\\
&\qquad
+48\U\tr\Lambda^0
\gint{\negquad\negquad
{\xy\rgvertex[{\theta_1}^{3}]{3}\loose1\endxy}\negquad\negquad}
\end{align*}
where
\begin{align*}
\overline{\psi}_3^{}(\theta_1,\theta_2,&\theta_3,\theta_4)=
\psi_3^{}(\theta_1,\theta_2,\theta_3)+\psi_3^{}(\theta_2,\theta_3,\theta_4)+
\psi_3^{}(\theta_3,\theta_4,\theta_1)+\\
&\phantom{mmmmmmmmmmmmmmmmmmmmm}
+\psi_3^{}(\theta_4,\theta_1,\theta_2)=\\
&=
\U\biggl(\frac{13}{2}{\theta_2}{\theta_1}^2-{\theta_2}{\theta_3}
{\theta_4}- {\theta_1}{\theta_2}{\theta_3}+7{\theta_3}^3+7
{\theta_1}^3
+7{\theta_2}^3+\\
&\qquad+7{\theta_4}^3+4{\theta_1}{\theta_3}^2
+4{\theta_3}{\theta_1}^2
-{\theta_3}{\theta_4}{\theta_1}
+\frac{19}{2}{\theta_1}{\theta_2}^2+\\
&\qquad
+\frac{19}{2}{\theta_2}{\theta_3}^2
+\frac{13}{2}{\theta_3}{\theta_2}^2
+\frac{13}{2}{\theta_1}{\theta_4}^2
+\frac{19}{2}{\theta_4}{\theta_1}^2-
\\&\qquad
-{\theta_4}{\theta_1}{\theta_2}
+4{\theta_4}{\theta_2}^2+\frac{19}{2}{\theta_3}{\theta_4}^2
+\frac{13}{2}{\theta_4}{\theta_3}^2+4{\theta_2}{\theta_4}^2
\biggr)
\end{align*}
and we have written the total contribution coming from terms
of total degree at least $4$ as a sum of terms of total
degree $3$.

The proof goes on along these lines: after a few more (completely straightforward,
but long and tedious) computations one ends up with an
explicit expression whose terms are all of total degree zero. These
terms cancel out and equation \prettyref{eq:ii} is proven. The interested reader can found this computation carried out in full detail in \cite{fiorenza-tesi}.

\section*{Appendix: Formal differential operators}\label{appendix}
In \prettyref{sec:witten} we showed how the second order derivative
$\partial^2Z(t_*)/\del t_0\del t_1$ appearing in equation
\prettyref{eq:ii} can be computed by considering Feynman diagrams with 
two distinguished holes.\par 
In this appendix we show how it can be computed by using the
main result from \cite{fiorenza-murri;matrix-integrals}. Let
$Z(s_*;t_*)$ be the partition function of the combinatorial
intersection numbers (see
\cite{kontsevich;intersection-theory;1992,arbarello-cornalba;dfiz}).
Then $Z(t_*)=Z(s^\circ_*;t_*)$, where $s_*^\circ=(0,1,0,0,\dots)$.
Also the combinatorial partition function $Z(s_*;t_*)$ is  related to
the asymptotic expansion of a matrix integral, namely
$$
Z(s_*;t_*)\biggr\rvert_{t_*(\Lambda)} =
\int_{\Hermitian[N]}\exp\left\{
    -\sqrt{-1}\sum_{j=0}^\infty(-1/2)^js_j\frac{\tr
      X^{2j+1}}{2j+1}
  \right\}\ud\mu_\Lambda(X),
$$
Moreover
$$
\frac{\del^{|n_*|} Z(s_*;t_*)} {{\del s_0}^{n_0} \cdots {\del
s_l}^{n_l}}
  \biggr\rvert_{s_*^\circ,t_*(\Lambda)}^{}\!\!=
\frac{n_0! \cdots n_l!}{\U^{|n_*|}(-2)^{\sum_j jn_j}}
\bigl\langle\!\bigl\langle{
    {{\sf v}_1}^{\coprod n_0} {\textstyle\coprod} \cdots
    {\textstyle\coprod} {{\sf v}_{2l+1}}^{\coprod n_l}
  }\bigr\rangle\!\bigr\rangle\!_{\Lambda}^{},
$$
where ${\sf v}_k$ denotes a $k$-valent special vertex and $\coprod$
denotes the disjoint union.
\par For any polyindex \(m_*=(m_0,m_1,\dots,m_l,0,0,\dots)\) set:
\begin{equation*}
\norm{m_*}_-:=\sum_{i=1}^\infty(2i-1)m_i,\quad
\norm{m_*}_+:=\sum_{i=0}^\infty(2i+1)m_i.
\end{equation*}
A formal triangular differential operator in the variables
$s_*$ is a formal series
  \begin{equation*}
    D(s_*,\del / \del{s_*}) = \sum_{\norm{n_*}_+ \leq \norm{m_*}_-}
    a_{m_*,n_*} s_*^{m_*} \frac{\del^{\abs{n_*}}}{\del s_*^{n_*}},
    \qquad a_{m_*, n_*} \in \setC,
  \end{equation*}
  of bounded degree in $s_*$.
 Formal triangular differential operators in the variables $s_*$ form a
(non commutative) algebra $\setC \langle\!\langle s_*,
\del / \del{s_*} \rangle\!\rangle$ which naturally acts on
$\setC[[t_*;s_*]]$. The main result from
\cite{fiorenza-murri;matrix-integrals} can stated as follows: there
exists formal triangular differential operators $D_k(s_*,s_*,\del /
\del{s_*})$ such that, for any $k_1,\dots, k_n$ in $\setN$,
$$
\frac{\del ^n Z(s_*;t_*)}{\del t_{k_1}\cdots \del t_{k_1}}=
D_{k_1}(s_*,\del / \del{s_*})\cdots
D_{k_n}(s_*,\del / \del{s_*})Z(s_*;t_*)
$$
As an immediate corollary, one gets
$$
\frac{\del ^2 Z(t_*)}{\del t_{1}\del t_{0}}
=
\bigl(D_{1}(s_*^\circ,\del / \del{s_*}) D_{0}(s_*,\del /
\del{s_*})Z(s_*;t_*)
\bigr)\biggr\rvert_{s_*^\circ}
$$
It is computed in
\cite{fiorenza-murri;matrix-integrals} that
\begin{equation*}
D_0(s_*,\del / \del{s_*})=\frac{{s_0}^2}{2}+\sum_{m=0}^{\infty}(2m+1)
s_{m+1}\frac{\del}{\del
s_{m}}
\end{equation*}
Similar computations give
\begin{equation*}
D_1(s_*^\circ,\del /
\del{s_*})=\frac{1}{2}
\frac{\del}{\del s_1}+\frac{1}{24}
\end{equation*}
We therefore find
\begin{align*}\label{eq:to-evaluate}
\frac{\del^2 Z(t_*)}{\del
      {t_1}\del t_0}\!\!&=\!
      \left.\left(\!\left(\frac{1}{2}
\frac{\del}{\del s_1}+\frac{1}{24}\right)\!\!\left(\frac{{s_0}^2}
{2}+\sum_{m=0}^{\infty}(2m+1)
s_{m+1}\frac{\del}{\del
s_{m}}\right)
Z(s_*;t_*)\right)\right|_{s_*^\circ}\\
&=\left.\left(\left(\frac{1}{2}\frac{\del^2}{\del s_0\del
s_1}+\frac{13}{24}\frac{\del}{\del
s_0}\right)Z(s_*;t_*)\right)\right|_{s_*^\circ}
\end{align*}
and we get the asymptotic expansion
$$
\frac{\del^2 Z(t_*)}{\del
      {t_1}\del t_0}\biggr\rvert_{t_*(\Lambda)}
=\frac{1}{4}\,\gint{\negquad\negquad\negquad{\xy\rgvertex[\!]{3}%
\loose1\endxy}\negquad\negquad\negquad{\xy\rgvertex[\!]{3}%
\loose1\loose2\loose3\endxy}}-\frac{13\U}{24}\,
\gint{\negquad\negquad\negquad{\xy\rgvertex[\!]{3}%
\loose1\endxy}\negquad\negquad\negquad}
$$
%%
%% Bibliografia
%%
%
%\bibliography{math.AG}
%\bibliographystyle{alpha}
\def\cprime{$'$} \def\cprime{$'$}

%\medskip
%
%\noindent{\it \textup{2000} Mathematics Subject Classification}: 81Q30
%(Primary); 14H81, 37K20 (Secondary).
%
%\medskip

%\noindent{\sc Dipartimento di Matematica ``Guido Castelnuovo'' ---
%Universit\`a degli Studi di Roma ``la Sapienza'' --- P.le Aldo Moro, 2 --
%00185 -- Roma, Italy}
%\par
%\noindent{\it E-mail address:} {\tt fiorenza@mat.uniroma1.it}
%

\end{document}